\documentclass[reqno, 12pt]{amsart}

\usepackage{amssymb}
\usepackage{latexsym}
\usepackage{epsfig}
\usepackage{epsf}
\usepackage{psfrag}
\usepackage{float}
\usepackage{subfigure}
\usepackage{a4wide}
\usepackage{hyperref}

\newtheorem{theorem}{Theorem}
\newtheorem{lemma}{Lemma}[section]

\newtheorem{proposition}{Proposition}[section]
\newtheorem{remark}{Remark}
\newtheorem{claim}{Claim}

\def\bR{{\mathbb R}}  
\def\bN{{\mathbb N}}  
\def\bP{{\mathbb P}}    
\def\bZ{{\mathbb Z}}  
\def\bQ{{\mathbb Q}} 
  
\def\bE{{\mathbb E}}
\def\bF{{\mathbb F}}
\def\bH{{\mathbb H}}
\def\bL{{\mathbb L}}

\def\1{{\mathbf 1}}
\def\0{{\mathbf 0}}
\def\A{{\mathbf A}}
\def\B{{\mathbf B}}
\def\C{{\mathbf C}}
\def\D{{\mathbf D}}

\def\S{{\mathbf S}}
\def\R{{\mathbf R}}
\def\Q{{\mathbf Q}}
\def\0{{\mathbf 0}}
\def\a{{\mathbf a}}
\def\e{{\mathbf e}}
\def\c{{\mathbf c}}

\def\v{{\mathbf v}}
\def\u{{\mathbf u}}
\def\x{{\mathbf x}}
\def\y{{\mathbf y}}
\def\z{{\mathbf z}}

\def\calF{{\mathcal F}}
\def\calN{{\mathcal N}}
\def\calT{{\mathcal T}}
\def\calU{{\mathcal U}}
\def\calV{{\mathcal V}}
\def\calE{{\mathcal E}}

\def\calC{{\mathcal C}}
\def\calB{{\mathcal B}}
\def\calD{{\mathcal D}}

\def\calP{{\mathcal P}}

\def\calR{{\mathcal R}}
\def\calS{{\mathcal S}}

\def\reff#1{(\ref{#1})}
\def\proofof #1{{\noindent \bf Proof of #1.}}
\def\endproof{$\square$ \vskip 2mm}

\begin{document}

\title[]{Multi-type shape theorems for FPP models}
\author{Leandro P. R. Pimentel}
\address{Institut de Math\'{e}matiques\\
\'{E}cole Polytechinique F\'{e}d\'{e}rale de Lausanne\\ 
CH-1015 Lausanne\\
Switzerland\\} 
\email{leandro.pimentel@epfl.ch}



\begin{abstract}
An Euclidean first-passage percolation (FPP) model describing the competing growth between $k$ different types of infection is considered. We focus on the long time behavior of this multi-type growth process and we derive multi-type shape results related to its morphology.
\end{abstract}

\maketitle

\section{Introduction}\label{int}

In standard planar first-passage percolation \cite{HW65} each pair $\x$ and $\y$ of nearest-neighbor of $\bZ^2$ has an edge connecting them and each edge is equipped with a non-negative random variable (passage time) which may be  interpreted as the time it takes for an infection to be transmitted from $\x$ to $\y$. We assume these random variables are i.i.d. with a continuous distribution $\bF$. The passage time $t(\gamma)$ for a nearest-neighbor path $\gamma$ is simply the sum of the passage times along the path. For $\x,\y\in\bZ^2$, the first-passage time from $\x$ to $\y$, which we denote $T(\x,\y)$, is the infimum of $t(\gamma)$ over all paths $\gamma$ from $\x$ to $\y$. For $t\geq 0$, let $\B(t)$ be the set of sites $\x$ reached from the origin $\0$ by time $t$, i.e. $T(\0,\x)\leq t$. One may think of sites in $\x\in \B(t)$ as infected and those in $\B(t)^c$ as healthy, and that at time $0$ the origin $\0$ is infected by some type of disease. The process $\big\{\B(t)\,:\,t\geq 0\big\}$ is then a model for the growth of an infection.

An interesting aspect of the evolution of the infection, namely the \emph{tree of infection}, is constructed as follows. First notice that, since the passage time distribution is continuous, for all $\x,\y\in\bZ^2$ there is (almost surely) an unique time-minimizing path (or geodesic) from $\x$ to $\y$, which we denote $\rho(\x,\y)$, such that $T(\x,\y)=t\big(\rho(\x,\y)\big)$. Thus $\rho(\x,\y)$ may be interpreted as the path through which the infection was transmitted from $\0$ to $\x$. With this picture in mind, the tree of infection $\Gamma$ is defined by the union of edges $\e\in\rho(\0,\x)$ over all $\x\in\bZ^2$. Newman \cite{n95} has shown that the number $K(\Gamma)$ of topological ends of $\Gamma$, i.e. the number of semi-infinite self-avoiding paths in $\Gamma$, is infinite provided an exponential moment condition on $\bF$ and a certain hypothesis concerning the uniformly bounded curvature of the asymptotic shape of $\B(t)$. In spite of the curvature hypothesis is plausible it has so far not been proved. 

\begin{figure}[htb]
\begin{center}
\includegraphics[width=0.5\textwidth]{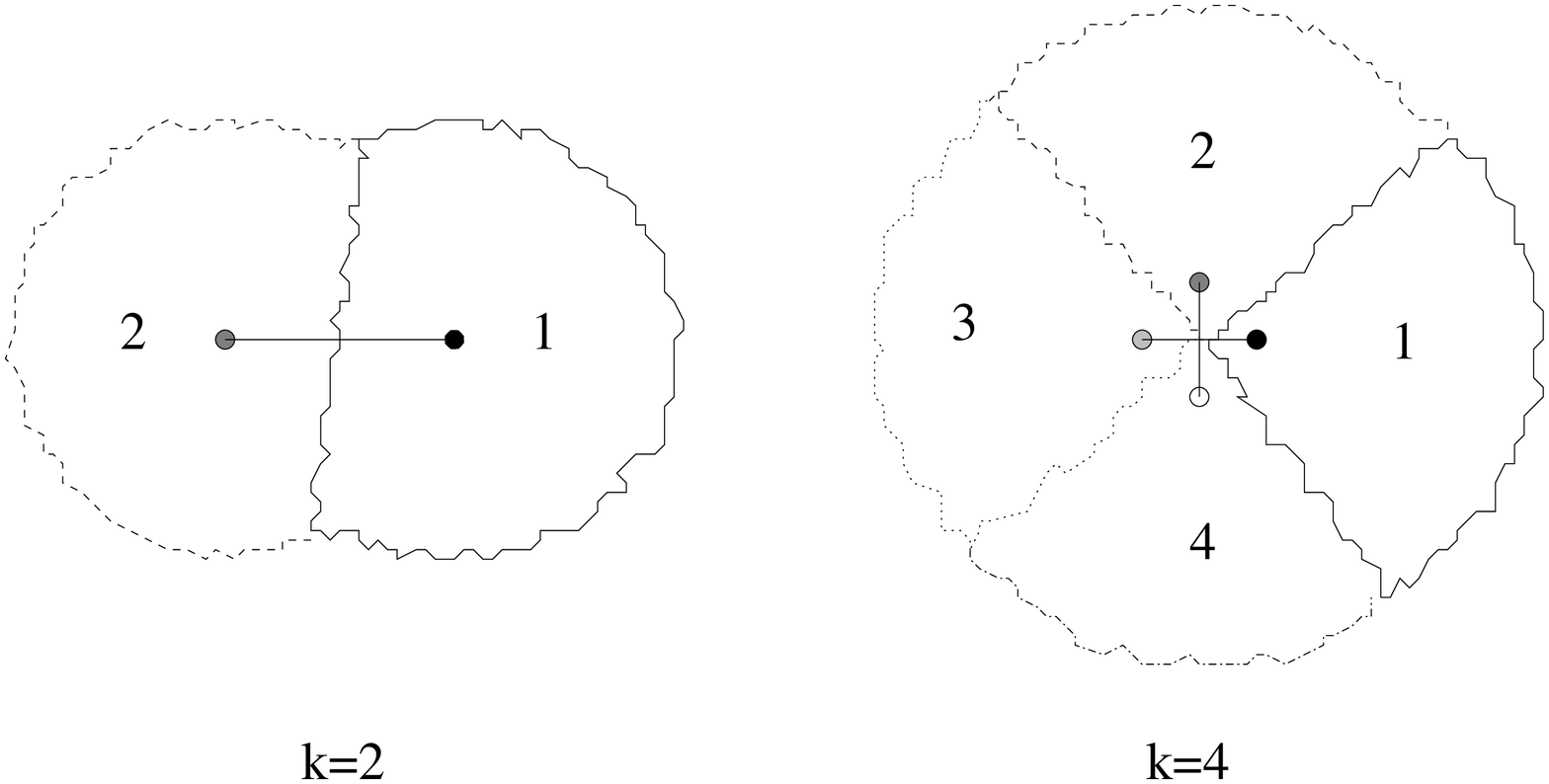}
\end{center}
\caption{Growth and Competition}\label{f2}
\end{figure}

In order to study the tree of infection, H\"agggstr\"om and Pemantle \cite{hp98,hp99} have introduced a multi-type growth model as follows. At time $0$ we start with $k$ different sites of $\bZ^2$, say $\x_1,\dots,\x_k$, each one representing a different type of infection. A site $\y\in\bZ^2$ is then infected at time $\min\big\{T(\x_1,\y),\dots,T(\x_k,\y)\big\}$ and it is acquired by the infection which first arrives there, i.e. by the unique type $j\in\{1,\dots,k\}$ such that $T(\x_j,\y)=\min\big\{T(\x_1,\y),\dots,T(\x_k,\y)\big\}$ (Figure \ref{f2})
. It may happens that at some early stage one of the types of infection completely surrounds another one, which then is prevented to grow indefinitely. If this does not occur, or equivalently, if all types of infection grow unboundedly, we say that $k$-coexistence occurs. 

Turning back to the question of topological ends of $\Gamma$, H\"aggstr\"om and Pemantle have noticed that if $k$-coexistence occurs with positive probability then $K(\Gamma)\geq k$ occurs with positive probability. They also  have shown that, if one considers an exponential passage time distribution then $2$-coexistence occurs with positive probability, and thus $K(\Gamma)\geq 2$ occurs with positive probability. Later Garet and Marchand \cite{gm04} and Hoffman \cite{H05} have extended this last result for stationary and ergodic FPP models on $\bZ^d$.

In this work we focus on the long time behavior of this multi-type growth model. However, differently from the above mentioned authors, we choose a first-passage percolation set-up on a random Delaunay triangulation \cite{VW90} whose spherical symmetry (isotropy) ensures that the asymptotic shape of the corresponding growth process is an euclidean ball. This choice allows us to prove various statements concerning minimizing paths, such as $\bP\big(K(\Gamma)=\infty\big)=1$, who could mostly only be conjectured by Newman in the standard model. In this setting, the main results we will prove are the following:
\begin{itemize}
\item If a type of infection survives then the region it conquers is (asymptotically) a cone with a random angle (Theorem \ref{t1}, Remark \ref{rand+strai});
\item If the $k$ initial sites form a regular polygon centered at the origin with radius $r$, then the probability that $k$ coexistence occurs tends to $1$ when $r$ tends to infinity. Moreover, for all $\epsilon>0$, the probability that for all $j\in\{1,\dots,k\}$ the region conquered by infection $j$ contains (asymptotically) the cone with axis through $\0$ and $\x_j$ and angle $\frac{\pi}{k}-\epsilon$ also tends to $1$ (Theorem \ref{t2}).
\end{itemize}

The main idea to prove our results is to explore the relation between this multi-type growth model and the asymptotic behavior of $T(\x,\y_n)-T(\0,\y_n)$ when $\y_n$ goes to infinity along a ray of angle $\alpha$ (Theorem \ref{tBuse-1} and Theorem \ref{tBuse-2}). We also study some roughening aspects of the one-dimensional boundary between the infections, namely the \emph{competition interface}, which were pointed out by physicists in numerical simulations \cite{dd91,sk95} (Remark \ref{rand+strai}). We note that analogous problems in the context of last-passage percolation and totally asymmetric exclusion processes were treated by Ferrari and Pimentel \cite{fp04-1} and Ferrari, Martin and Pimentel \cite{fp04-2}. Deijfen, H\"aggstr\"om and Bagley \cite{dhb03} have also considered isotropic multi-type growth models in $\bR^d$ where the growth is driven by outbursts in the infected region. 

\subsection{Multi-type growth process} 
\begin{figure}[htb]
\begin{center}
\includegraphics[width=0.3\textwidth]{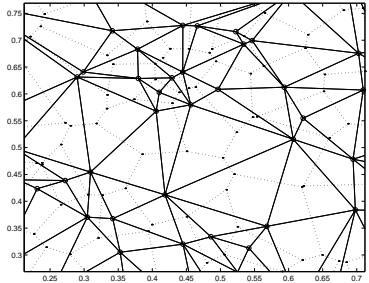}
\end{center}
\caption{The Delaunay triangulation and the Voronoi tessellation}\label{f1} 
\end{figure}

Consider the random graph $\calD:=(\calD_v,\calD_e)$, named the \emph{Delaunay triangulation}, constructed as follows. The vertex set $\calD_v\subseteq\bR^2$ is the set of points realized in a two-dimensional homogeneous Poisson point process with intensity $1$. To each vertex $\v$ corresponds an open and bounded polygonal region $\C_\v$ (the Voronoi tile at $\v$) consisting of the set of
points of $\bR^2$ which are closer to $\v$ than to any other $\v'\in\calD_v$. The
edge set $\calD_e$ consists of non oriented pairs $(\v,\v')$ such that $\C_\v$ and
$\C_{\v'}$ share a one-dimensional edge (Figure \ref{f1}). One can see that (with probability one) each Voronoi tile is a convex and bounded polygon, and the graph $\calD :=(\calD_v,\calD_e)$ is a triangulation of the plane  \cite{M91}.
The \emph{Voronoi tessellation} $\calV:=(\calV_v,\calV_e)$ is defined by
taking the vertex set $\calV_v$ equal to the set of vertices of the Voronoi tiles and
the edge set $\calV_e$ equals to the set of edges of the Voronoi tiles. 

Each edge $\e\in\calD_e$ is
independently assigned a nonnegative random variable $\tau_\e$ from a common
distribution $\bF$ (the passage time distribution) that is independent of the
Poisson process $\calD_v$. We assume throughout that $\bF$ is continuous and that
\begin{equation}\label{a1}
\int e^{ax}\bF(dx)<\infty\,\mbox{ for some }\,a\in(0,\infty)\,.
\end{equation}
We denote by $(\Omega,\calF,\bP)$ our underline probability space, i.e. from each realization $\omega\in\Omega$ one can determine the Poisson point process as well the passage time configuration. This model inherits the euclidean  (translation and rotational) invariance of the Poisson point process.

The passage time $t(\gamma)$ of a path $\gamma$ in $\calD$ is the sum of the
passage times of the edges in $\gamma$:
\[
t(\gamma):=\sum_{\e\in\gamma}\tau_\e \,.
\]
The first-passage time between
two vertices $\v$ and $\v'$ in $\calD_v$ is defined by
\[
T(\v,\v'):=\inf\big\{t(\gamma);\,\gamma\in\calC(\v,\v')\big\} ,
\]
where $\calC(\v,\v')$ the set of all paths connecting $\v$ to $\v'$. We extend the first-passage time $T$ to $\x,\y\in\bR^2$ by setting $T(\x,\y):=T\big(\v(\x),\v(\y)\big)$, where $\v(\x)$ is the almost sure unique vertex $\v\in\calP$ with $\x\in\C_{\v}$. We say that $\rho(\v,\v')\in\calC(\v,\v')$ is a geodesic between $\v$ and $\v'$ if $t\big(\rho(\v,\v')\big)=T(\v,\v')$. For each $\x,\y\in\bR^2$ we denote $\rho(\x,\y):=\rho\big(\v(\x),\v(\y)\big)$. One can see that if $\bF$ is a continuous function then, almost surely, for all $\v,\v'\in\calD_v$ there exists a unique geodesic $\rho(\v,\v')$ \cite{p-105}. A self-avoiding and semi-infinite path $\rho=(\v_1,\v_2,\dots)$ in $\calD$ is called a semi-infinite geodesic if for all $\v_j,\v_k\in\rho$, the path $(\v_j,\v_{j+1},...,\v_k)$ is the unique geodesic connecting $\v_j$ to $\v_k$.

Given $k$ different points $\x_1,...,\x_k\in\bR^2$, the initial configuration of seeds, we define the multi-type growth process $\big\{(\B_{\x_1}(t),...,\B_{\x_k}(t))\,:\,t\geq 0\big\}$ by 
\[
 \B_{\x_j}(t):=\big\{\x\in\bR^2\,:\,\x\in c(\C_\v)\,\mbox{ for some }\,\v\in\calB_{\x_j}(t)\big\}\,,
\]
where
\[
 \calB_{\x_j}(t):=\big\{\v\in\calD_v\,:\,T(\x_j,\v)\leq t\,\mbox{ and }\,\min_{l=1,...,k}\{T(\x_l,\v)\}=T(\x_j,\v)\big\}\,,
\]
and $c(\C_\v)$ denotes the closure of the tile $\C_\v$. If there exists $j<l$ such that
$\v_{\x_j}=\v_{\x_l}$ then we set $\B_{\x_j}$ as before and $\B_{\x_l}(t)=\emptyset$. 

When $k=1$ then we have a single
growth process $\B_\x(t)$ which represents the set of points reached by
time $t$ from the initial seed $\x$. For a continuous distribution $\bF$ satisfying \reff{a1} the following shape theorem \cite{p-105,VW92} holds: there exists a constant  $\mu(\bF)\in(0,\infty)$, namely the time constant, such that for all $\epsilon>0$
\[
\bP\big((1-\epsilon )t \D(1/\mu)\subseteq \B_\0(t)\subseteq (1+\epsilon )t
\D(1/\mu)\mbox{ eventually }\big)=1 \,,
\]
where $\D(r):=\{\x\in\bR^2 \,:\, |\x|\leq r\}$ and $\0:=(0,0)$. 

When $k\geq2$ the process $\big\{(\B_{\x_1}(t),...,\B_{\x_k}(t))\,:\,t\geq 0\big\}$ is a model for competing growth on the plane where each point $\x\in\bR^2$ is acquired by the specie $j\in\{1,\dots,k\}$ which first arrives there. The competition interface $\psi$ is the one-dimensional boundary between the species when $t=\infty$. This interface can be seen as a finite union of polygonal curves determined by edges in $\calV$ (the Voronoi tessellation) which are shared by tiles in different species. A branch of the competition interface is a self-avoiding path $\varphi=(\x_n)_{n\geq 1}$ in $\calV$ such that $\{\x_n\,:\,n\geq 1\}\subseteq\psi$.   

For each $\alpha\in[0,2\pi)$ we say that a self-avoiding path $(\x_{n})_{n\geq 1}$, with vertices in $\bR^2$ and such that $|\x_n|\to\infty$ when $n\to\infty$, is a $\alpha$-path if
\[
\lim_{n\to\infty}\frac{\x_{n}}{|\x_{n}|}=e^{i\alpha}:=(\cos\alpha,\sin\alpha)\,. 
\]
In this case we also say that $(\x_n)_\bN$ has the asymptotic orientation $e^{i\alpha}$. This is equivalent to  
\[
\lim_{n\to\infty}ang(\x_{n},e^{i\alpha})=0\,,
\]
where $ang(\x,\y)$ denotes the angle in $[0,\pi]$ between the points $\x,\y\in\bR^2$. Thus, a sufficient condition for a path $(\x_{n})_{n\geq 1}$ to be a $\alpha$-path for some $\alpha\in[0,2\pi)$ is, for some fixed $\delta\in(0,1)$ and some constant $c>0$, for sufficiently large $n$
\[
 ang(\x_{n},\x_m)\leq |\x_n|^{-\delta}\,\mbox{ whenever }m > n\,,
\]
which is the so called $\delta$-straightness property for semi-infinite paths introduced by Newman \cite{n95}.


\begin{theorem}\label{t1}
For $k\geq 2$ let $\Omega_k$ be the event that, for the competing growth model with $k$-different species, there exists a finite subset $\Theta:=\{\theta_1,...,\theta_m\}$ of $[0,2\pi)$ such that every branch $\varphi$ of the competition interface is a $\theta(\varphi)$-path for some $\theta\in\Theta$. Under \reff{a1}, $\bP\big(\Omega_k\big)=1$.
\end{theorem} 

\begin{remark}\label{rand+strai}
In Section \ref{geo} (part \ref{pr-rand+strai}) we will give a sketch of the proof that for all  $\alpha\in[0,2\pi)$ 
\[
\bP\big(\theta=\alpha\mbox{ for some }\theta\in\Theta\big)=0\,,
\] 
and that if $\xi\in(3/4,1)$ then, almost surely, for all branch $\varphi=(\x_{n})_{n\geq 1}$ of the competition interface there is a constant $c>0$ such that 
\[
ang(\x_{n},e^{i\theta(\varphi)})\leq c|\x_n|^{\xi-1}\mbox{ eventually }\,.
\]
\end{remark}

Let $\x_1(r)=(0,r),\dots,\x_k(r)$ be the vertices of a regular polygon with $k$ sides and radius $r$. For each $j=1,...,k$ define
the projection of the random set $\B_{j}^r:=\B_{\x_j(r)}(\infty)$ onto $\S^1$, the set of unit vectors $|\x|=1$, by 
\[
\S_{j,r}:=\{\x = e^{i\alpha}\in \S^{1}\,:\,\bL_{s\x}(\alpha)\subseteq \B_{j}^{r}\,\mbox{ for some }\,s>0\}\,,
\]
where $\bL_\x(\alpha)$ denotes the line starting from $\x$ and with direction $e^{i\alpha}$. For each $\epsilon\in (0,\pi/k)$ and $j\in\{1,\dots,k\}$  define 
\[
\S_{j}(\epsilon):=\{\x\in\S^1\,:\, ang(\x,\x_j(r))\leq \frac{\pi}{k}-\epsilon\}\,. 
\]

\begin{theorem}\label{t2}
Let $k\geq 2$. Under \reff{a1}, for all $\epsilon>0$
\[ 
 \lim_{n\to\infty}\bP\big(\S_{j}(\epsilon)\subseteq \S_{j,n}\,\emph{ for all $j=1,\dots,k$ }\big)=1\,.
\]
\end{theorem}


\subsection{Busemann type asymptotics and the competition interface}
To illustrate the approach we follow in this work to study the competition interface assume that $k=2$. Consider the line $\bL_\0(\alpha)$ starting from the origin $\0$ and with direction $e^{i\alpha}$. Then we have three possibilities: i) either it intersets the competition interface infinitely many times; ii) or it is eventually contained in  $\B_{\x_1}(\infty)$; iii) or it is eventually contained in $\B_{\x_2}(\infty)$. Notice that the former  implies 
\[
\liminf_{s\to\infty}\big(T(\x_{1},se^{i\alpha})-T(\x_{2} ,se^{i\alpha})\big)\leq
0\leq\limsup_{s\to\infty}\big(T(\x_{1},se^{i\alpha})-T(\x_{2},se^{i\alpha})\big)\,,
\]
while the second implies 
\[
\limsup_{s\to\infty}\big(T(\x_{1},se^{i\alpha})-T(\x_{2} ,se^{i\alpha})\big)\leq 0\,,
\]
and the third implies
\[
0\leq \liminf_{s\to\infty}\big(T(\x_{1},se^{i\alpha})-T(\x_{2} ,se^{i\alpha})\big)\,.
\]
It turns out that the above expressions resemble Busemann type asymptotics for $T$ (see Ballmann \cite{b95} for more details on this subject). Newman \cite{n95,ln96} has shown for the lattice model that, under suitable assumptions on the curvature of the limit shape, $T(\x_1,\y_n)-T(\x_2,\y_n)$ attains eventually a nonzero value $H^{\alpha}(\x_1,\x_2)$, called the Busemann function. By following his method, and by taking profit of the isotropy in our model, we will show\footnote{We also refer to \cite{hn97}, where an analog result is proved in an Euclidean first-passage percolation set-up.}:

\begin{theorem}\label{tBuse-1}
For $\alpha\in[0,2\pi)$ let $\Omega_0(\alpha)$ be the event that for all $\v,\bar{\v}\in\calD_v$, there exists $H^{\alpha}(\v,\bar{\v})$, nonzero for $\v\neq\bar{\v}$, such that 
\begin{equation}\label{eBuse-1}
\lim_{{|\x|\to\infty}\atop{\x/|\x|\to e^{i\alpha}}}\big(T(\v,\x)-T(\bar{\v},\x)\big)=H^{\alpha}(\v,\bar{\v})\,.
\end{equation}
Under \reff{a1}, $\bP\big(\Omega_0(\alpha)\big)=1$.
\end{theorem}


For $\x ,\y\in\bR^2$ we set $H^\alpha(\x,\y):=H^\alpha\big(\v(\x),\v(\y)\big)$. It was conjectured by Howard and Newman \cite{hn01} that
\[
\lim_{n\to\infty}\frac{H^{\alpha}(n\vec{e}_1,\0)}{n}=-\mu(\bF)\cos\alpha\,, 
\]
where $\vec{e}_1:=(1,0)$. This observation is related to the asymptotic behavior of our multi-type growth model and the key result to show Theorem \ref{t2} is the following theorem, which is a small step towards the above conjecture.

\begin{theorem}\label{tBuse-2}
For $\alpha\in[0,\pi/2)$ let $\Omega_1(\alpha)\subseteq\Omega_0(\alpha)$ be the event that 
\begin{equation}\label{eBuse-2}
-\mu(\bF)\leq \lim\inf_{n\to\infty}\frac{H^{\alpha}(n\vec{e}_{1},\0)}{n}\leq \limsup_{n\to\infty}\frac{H^{\alpha}(n\vec{e}_{1},\0)}{n}\leq -\mu(\bF)\frac{\cos\alpha}{1+\sin\alpha}\,.
\end{equation}
Under \reff{a1}, $\bP\big(\Omega_1(\alpha)\big)=1$. In particular, with probability one,
\[
\lim_{n\to\infty}\frac{H^{0}(n\vec{e}_{1},\0)}{n}=-\mu(\bF)\,.
\]
\end{theorem}

\subsection*{Overview} In Section \ref{multi} we will deduce Theorem \ref{t1} and Theorem \ref{t2} from Theorem \ref{tBuse-1} and Theorem \ref{tBuse-2}. In Section \ref{pre} we will start by defining the probability space where our model takes place and we will show a modification lemma that will play an important rule in the study of coalescence of semi-infinite geodesics. After that we will study some geometrical aspects of Voronoi tilings. We note that in the Delaunay triangulation context some technical difficulties are imposed by its long range dependence. Some of them will be avoided by making references to results of a previous work of the author \cite{p-105,p-205}. In the third part we will recall some geometrical lemmas concerning the $\delta$-straightness of semi-infinite paths. Finally, in Section \ref{geo} will study existence and coalescence of semi-infinite geodesics to show Theorem \ref{tBuse-1} and Theorem \ref{tBuse-2}. It  will largely parallel the analog study develop by Newman et al \cite{hn01, ln96, n95, np95} in the lattice and in the Euclidean FPP models.

\section{Proof of the multi-type shape theorems}\label{multi}

\proofof{Theorem \ref{t1}}  For each $j=1,...,k$, let $\S_{j}$ denote the set of unit vectors $e^{i\beta}$
such that $\bL_{se^{i\beta}}(\beta)\subseteq \B_{\x_j}(\infty)$ for some $s>0$ and let 
\[
\S_{0}:=(\cup_{j=1}^{l}\S_{j})^{c}\,. 
\]
Let 
\[
\D_n:=\{e^{i\beta}\,:\,\beta=2k\pi/2^{n}\mbox{ for some }1\leq k\leq 2^{n}\}\, 
\]
and $\D:=\cup_{n\geq 1}\D_n$. Consider the event $\cap_{\alpha\in\D}\Omega_0(\alpha)$ that for all $\alpha\in\D$ and $\v,\bar{\v}\in\calD_v$  there exists $H^{\alpha}(\v,\bar{\v})$, nonzero for $\v\neq\bar{\v}$, such that 
\[
\lim_{{|\x|\to\infty}\atop{\x/|\x|\to e^{i\alpha}}}\big(T(\v,\x)-T(\bar{\v},\x)\big)=H^{\alpha}(\v,\bar{\v})\,.
\]
By Theorem \ref{tBuse-1}, $\bP\big(\cap_{\alpha\in\D}\Omega_0(\alpha)\big)=1$. 

We claim that, on this event, every branch of the competition interface is an $\theta$-path for some $\theta\in[0,2\pi)$. To see this, notice that if $e^{i\alpha}\in \S_{0}$ then for some $j_1 \ne j_2$, $\bL_{\0}(\alpha)$ intersects infinitely many times the region
$\B_{\x_{j_1}}(\infty)$ and the region $\B_{\x_{j_2}}$. Thus
\[
\liminf_{s\to\infty}\big(T(\x_{j_1},se^{i\alpha})-T(\x_{j_2} ,se^{i\alpha})\big)\leq
0\leq\limsup_{s\to\infty}\big(T(\x_{j_1},se^{i\alpha})-T(\x_{j_2},se^{i\alpha})\big)\,,
\]
which implies that $\D\cap\S_0=\emptyset$. Let $\C_k^n$ be the cone consisting of points $re^{i\beta}$ such that $r>0$ and $\beta\in(2\pi k/2^n,2\pi(k+1)/2^n)$. Now, if $\D\cap\S_0=\emptyset$ and $e^{i\beta}\in\D$ then every branch $\varphi$ of the competition interface can not intersect infinitely many times the line $\bL_{\0}(\beta)$. So,  for each branch $\varphi$ of the competition interface we can find a sequence of cones  $(\C^{n}_{k_{n}})_{n\geq1}$, with $n\to\infty$ and $\C^{n+1}_{k_{n+1}}\subseteq \C^{n}_{k_{n}}$, such that $\varphi$ is eventually contained in $C^{n}_{k_{n}}$. This implies that $\varphi$ must be a $\theta$-path for some $\theta\in[0,2\pi)$. \begin{flushright}\endproof\end{flushright}

\proofof{Theorem \ref{t2}} Since
\[
 \bP\big(\S_{j}(\epsilon)\subseteq \S^{r}_{j}\big)=\bP\big(\S_{1}(\epsilon)\subseteq \S^{r}_{1}\big)
\]
for all $j=1,...,k$, we only need to prove that
\begin{equation}\label{e1kshape}
 \lim_{r\to\infty}\bP\big(\S_{1}(\epsilon)\subseteq \S^{r}_{1}\big)=1\,.
\end{equation}
To do so, for each $j=1,\dots,k$ let $\alpha^k_j:=\pi(j-1)/k$, $\vec{e}^k_j:=e^{2i\alpha^k_j}$ and  $A_{r}:=\cap_{j=1}^{k} A^{j}_{r}$, where
\[
A_{r}^{j}:= \cap_{l\ne j}\big[H^{\alpha^k_j}(r\vec{e}^k_l,r\vec{e}^k_j)>0\big]\,.
\]
Let $\alpha^{+\epsilon}_{k}:=\frac{\pi}{k}-\epsilon$ and $\alpha^{-\epsilon}_{k}:= (2\pi-\frac{\pi}{k})+\epsilon$ and set
\[
B_{r}(\epsilon):=\cap_{j=2,\dots,k}\big[H^{\alpha^{+\epsilon}_{k}}(r\vec{e}^k_j,r\vec{e}^k_1)>0\mbox{ and }H^{\alpha^{-\epsilon}_{k}}(r\vec{e}^k_j,r\vec{e}^k_1)>0\big]\,. 
\]
By Theorem \ref{tBuse-2}
\begin{equation}\label{e2kshape}
\lim_{r\to\infty}\bP\big(A_{r}\cap B_r(\epsilon)\big)=1\,.
\end{equation}
The connectivity of the regions $\B_j^r$ yields that, on $A_r\cap B_r(\epsilon)$, $\S_1(\epsilon)\subseteq \S_1^r$. Together with \reff{e2kshape}, this yields \reff{e1kshape} and the proof of Theorem \ref{t2} is complete. \begin{flushright}\endproof\end{flushright}

\section{Auxiliary results}\label{pre}

\subsection{The probability space}\label{pre-ps}
During the subsequent proofs we will consider the following construction of $(\Omega,\mathcal{F},\bP)$, the underline probability space of our FPP model. Let $\u_{0}=(0,0),\u_{2},\dots$ be a ordering of
$\bZ^{2}$ and for each $k\geq 1$ let
\[
\B_{k}:=\u_{k}+[-1/2,1/2]^{2}\,.
\]
Consider 
\[
\calN=\{N_{k}\,:\, k\geq 1\},
\]
a collection of i.i.d. Poisson random variables with intensity $1$;
\[
\calU_k=\{U_{k,l}\,:\, l\geq 1\}, 
\]
a collection of independent random points in the plane so that $U_{k,l}$ has an
uniform distribution in the square box $B_{k}$;
\[
\calT_k=\{\tau_{k,l}^{m,n}\,:\, l\geq 1,m\geq k,n\geq 1\mbox{ and }n>l\mbox{
  whenever }k=m\}, 
\]
a collection of i.i.d. non negative random variables with common distribution
$\bF$ (the passage time distribution). We also impose that all these
collections are independent of each other.

To determine the vertex set $\calD_v=\calP$,
at each square box $\B_{k}$ we put $N_{k}$ points given by
$U_{k,1},...,U_{k,N_{k}}$. This procedure determines a Poisson point process
$\calP$ from the collections $\calN$ and $\calU_k$ with $k\geq 1$. Given $\e\in\calD_e$ we know that there exist an unique pair $(U_{k,l},U_{m,n})$, where
either $m>k$ or $m=k$ and $n>l$, so that $\e=(U_{k,l},U_{m,n})$. Set $\tau_{e}=\tau_{k,l}^{m,n}$.

For each $k\geq 1$ denote by $(\Omega^{k},\mathcal{F}^{k},\bP^{k})$ the probability induced by the random variable $N_k$ and the collections  $\calU_k$, $\calT_k$. The  probability space $(\Omega,\calF,\bP)$ is defined to be the product space of $(\Omega^{k},\calF^k,\bP^{k})$ over $k\geq 1$.

An important step in the construction of the Busemann function is the proof of the coalescence behavior of semi-infinite geodesics with the same asymptotic orientation. In this proof, we will use the following modification lemma. Let $\mathnormal{K}$ be the collection of all finite sequences 
\begin{equation}\label{E:prescription}
 I=\big((k_{j},l_{j},m_{j},n_{j})\big)_{j=1,...,q}\in(\bN^{4})^{q}
\end{equation}
where $q\geq 1$, $(k_{j},l_{j},m_{j},n_{j})\ne(k_{i},l_{i},m_{i},n_{i})$ for $j\ne i$,
$k_{1}\leq...\leq k_{q}$, and either $ k_{j}< m_{j}$ or $l_{j}<n_{j}$. To each
$I\in\mathnormal{K}$ corresponds a random vector $(\tau_{k_{j},l_{j}}^{m_{j},n_{j}})_{j=1,...,q}$. We denote  $(\Omega_{I},\calF_{I},\bP_{I})$ the
probability space induced by this random vector. Let
\[
\hat{\Omega}_{I}:=\{\hat{\omega}_{I}\,:\,\exists\,\omega_{I}\in\Omega_{I}\mbox{ with }(\hat{\omega}_{I},\omega_{I})\in\Omega\}
\]
and denote by $\hat{\bP}_I$ the probability law $\bP$ restricted to this subset. For each $I\subseteq\mathnormal{K}$, $A\subseteq\Omega$ and
$\omega_{1}\in\hat{\Omega}_{I}$ define
\[
A_{I,\omega_{1}}:=\{\omega_{2}\in\Omega_{I}\,:\,\omega=(\omega_{1},\omega_{2})\in A\}\,. 
\]
Let $\{R_{I}\,:\,I\in\mathnormal{K}\}$ be a family of events $R_{I}\in\calF_{I}$ such that $\bP_{I}(R_{I})>0$ for all $I$. Then define the map on $\calF$ by
\[
 \Phi_{I}(A):=\{\omega_{1}\in\hat{\Omega}_{I}\,:\,\bP_{I}(A_{I,\omega_{1}})>0\}\times R_{I}.
\]
Suppose that $W(\omega)$ is a random element of $\mathnormal{K}$, which may be interpreted as the set of indexes (edges) whose passage time value will be modified. For $A\subseteq\Omega$, let
\[
 \tilde{\Phi}(A):=\cup_{I\in\mathnormal{K}}\big[\{\omega_{1}\in\hat{\Omega}_{I}\,:\,A(I)_{I,\omega_{1}}\ne\emptyset\}\times R_{I}\big],
\]
where $A(I):=A\cap [W=I]$. 
\begin{lemma}\label{lmodif} For each $A\in\calF$, $\tilde{\Phi}(A)$ contains $\Phi(A)\in\calF$ defined as the following union 
\[ 
 \Phi(A):=\cup_{I\in\mathnormal{K}}\Phi_{I}\big(A(I)\big).
\]
Furthermore, if $\bP(A)>0$ then $\bP\big(\Phi(A)\big)>0$.
\end{lemma} 

\proofof{Lemma \ref{lmodif}} If $\bP_{I}\big(A(I)_{I,w_{1}}\big)>0$ then $A(I)_{I,w_{1}}\ne\emptyset$ and so  $\Phi(A)\subseteq\tilde{\Phi(A)}$. Since $\mathnormal{K}$ is countable and $A=\cup_{I\in\mathnormal{K}}A(I)$, if $\bP(A)>0$ then there exists $I\in\mathnormal{K}$ such that $\bP\big(A(I)\big)>0$. For this $I$, by Fubini's theorem
\begin{equation}\label{emodi2} 
 0<\bP\big(A(I)\big)=\int_{\hat{\Omega}_{I}}\bP_{I}\big(A(I)_{I,w_{1}}\big)\hat{\bP}_{I}(dw_{1}).
\end{equation}
Let 
\[
\hat{A}_{I}:=\{ w_{1}\,:\,\bP_{I}\big(A(I)_{I,w_{1}}\big)>0\}\,. 
\]
By \reff{emodi2}, $\hat{\bP}_{I}\big(\hat{A}_{I}\big)>0$. According to the definition of
$\Phi_{I}$, 
\[
\bP\Big(\Phi_{I}\big(A(I)\big)\Big)=\hat{\bP}_{I}(\hat{A}_{I})\bP_{I}(R_{I})>0\,. 
\]
Since $\Phi_{I}\big(A(I)\big)\subseteq\Phi(A)$, we conclude that $\bP\big(\Phi(A)\big)>0$. \begin{flushright}\endproof\end{flushright}

\subsection{Some geometrical aspects of Delaunay triangulations}\label{pre-geom}
In this part we are going to study some geometrical aspects of Delaunay triangulations.  Let $\x,\y\in\bR^2$ and construct a path $\gamma(\x,\y):=(\v_1,...,\v_{k})$ in $\calD$ connecting $\v(\x)$ to $\v(\y)$ as follows: set $\v_1:=\v(\x)$; if $\v_1\neq\v(\y)$ let $\v_{2}$ be the (almost-surely) unique nearest neighbor of $\v_1$ such that the edge of $\C_{\v_1}$ that is perpendicular to the line segment $[\v_1,\v_2]$ cross $[\x,\y]$; given $\v_l$ with $l\geq 1$, if $\v_{l}\neq\v(\y)$ then we set $\v_{l+1}$ to be the (almost-surely) unique nearest neighbor of $\v_l$, different from $\v_{l-1}$, such that the edge of $\C_{\v_l}$ that is perpendicular to  $[\v_l,\v_{l+1}]$ cross $[\x,\y]$; otherwise we set $k:=l$ and the construction is finished. We denote  $|\gamma(\0,n\vec{e}_1)|$ the number of edges in $\gamma(\0,n\vec{e}_1)$.

For $\z\in\bR^2$ and $L>0$ let 
\[
\B_\z^{L}:=L\z+[-L/2,L/2]\,
\]
For $n>0$ consider the set $\calE_n$ composed of edges $(\v,\bar{\v})\in\calD_e$ with $\C_\v\cap\B_\z^{1}\neq\emptyset$ or $\C_{\bar\v}\cap\B_\z^{1}\neq\emptyset$ for some $\z\in[\0,n\vec{e}_1]$. We denote $|\calE_n|$ the number of edges in $\calE_n$.
\begin{lemma}\label{lgraph2}
There exists constants $z_j,c_j>0$ such that for all $n\geq 1$,
\begin{equation}\label{egraph2*}
\bP\big(|\gamma(\0,n\vec{e}_1)|\geq z n\big)\leq e^{-c_1 zn}\,\mbox{ whenever }z\geq z_0\,,
\end{equation}
and
\begin{equation}\label{egraph3*}
\bP\big(|\calE_n|\geq z n\big)\leq e^{-c_2 zn}\,\mbox{ whenever }z\geq z_1\,.
\end{equation}
\end{lemma}

The proof of this lemma is performed through renormalization ideas developed in \cite{p-205}. To avoid some repetitions we give a sketch of the proof and leave the details to the reader, which can be filled by following the arguments in the proof of Proposition 2.2 in \cite{p-205} (which is exactly the proof in \reff{egraph2*}).

\begin{figure}[htb]
\begin{center}
\includegraphics[width=0.2\textwidth]{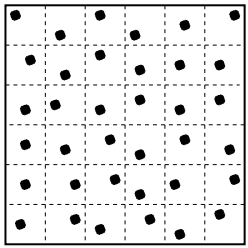}
\end{center}
\caption{Renormalization: a full box}\label{ffull}
\end{figure}

\proofof{Lemma \ref{lgraph2}} For $\z\in\bZ^2$ and $L>0$ divide a square box $\B_\z^{L}$ into thirty-six sub boxes of the same length, say $\B_1,\dots,\B_{36}$. We stipulate $\B$ is a full box if all those thirty-six sub boxes have at least one Poissonian point (Figure \ref{ffull}).

We say that $\Lambda:=(\B_{\z_1}^{L},\dots,\B_{\z_k}^{L})$ is a circuit of boxes if $(\z_1,\dots,\z_k)$ is a circuit in $\bZ^2$ (in the usual sense). Let $\lambda$ be the closed polygonal path composed by the line segments connecting $L\z_j$ to $L\z_{j+1}$, where $j=1,\dots,k-1$, together with $[\z_k,\z_1]$. To each circuit $\Lambda$ we associate two subsets of the plane: $\Lambda^{in}$ denotes the interior of the bounded component of $\bR^2\backslash\cup_{j=1}^k\B_{\z_j}^{L}$ while $\lambda^{in}$ will denote the interior of the bounded component of $\bR^2\backslash\lambda$. Now, assume that $\Lambda:=(\B_{\z_j}^{L})_{j=1}^k$ is a circuit composed by full boxes. By Lemma 2.1 in \cite{p-205}, we have the following geometrical property: if $\C_{\v}\cap \Lambda^{in}\neq\emptyset$ then $\C_{\v}\subseteq \lambda^{in}$. One important consequence of this is that the set of vertices used by $\gamma(\0,n\vec{e}_1)$ or by $\calE_n$ are both contained in the region $\R_n$ limited by the smallest circuit of full boxes surrounding the line segment $[\0,n\vec{e}_1]$. Therefore to show Lemma \ref{lgraph2} is enough to prove the analog decay for the number of Poissonian points in $\R_n$ \footnote{ Recall that, by the Euler formula, the number of edges and vertices in a triangulation have the same order.}.

Notice that, since each box is full independently of each other and the probability that it occurs goes to $1$ when $L$ goes to infinity, for a fixed large $L_0>0$, the probability that $\R_n$ contains more than $zn$ boxes decays as $e^{-czn}$ (see for instance Grimmet \cite{G99}).

Now, the number of points in $\R_n$, say $R_n$, is the sum of independent Poisonian random variables. This is less or equal to $M_m$, the maximum of the number of points in $\R$ over all connected regions $\R$ intersecting at most $m$ boxes the origin $\0$. Thus, on the event that $\R_n$ contains less than $zn$ boxes, we have $R_n\leq M_{zn}$. On the other hand,  $M_m$ can be seen as a Greedy lattice animal model and for such a model we can also show, for large $\bar{c}>0$, that the probability that $M_{m}\geq \bar{c}m$ decays as $e^{-cm}$ (Lemma 2.3 of \cite{p-205}).

By cooking together the arguments in these two last paragraphs one obtains that the probability that the number of points in $\R_n$ is greater than $zn$ also decays as $e^{-czn}$, for some constant $c>0$ and sufficiently large $z$.
\begin{flushright}\endproof\end{flushright}

Let $T_\calD$ denote the graph
metric on $\calD$, i.e. for $\v,\bar{\v}\in\calD_v$,
$T_\calD(\v,\bar{\v})$ is the minimum number of edges that one path should pass to go from
$\v$ to $\bar{\v}$. Notice that $T_\calD(\v,\bar{\v})$ is the first-passage time between $\v$ and $\bar{\v}$ if one associates to each edge $\e$ the passage time value $1$. For each $\A,\B\subseteq\bR^2$ we set $T_\calD(\A,\B)$ to be the minimum of $T_\calD(\v,\tilde{\v})$ over all pairs $\v$ and $\tilde{\v}$ such that $\C_\v\cap\A\neq\emptyset$ and $\C_{\tilde{\v}}\cap\B\neq\emptyset$. By the shape theorem, we have:
\begin{lemma}\label{lgraph1}  
There exists $\nu\in(0,\infty)$ such that almost surely
\[
\lim_{n\to\infty}\frac{T_\calD(\A,ne^{i\alpha}+\B)}{n}=\nu\,.
\]
\end{lemma}

We notice that $\nu$ does not depend either on $\A$ and $\B$ or on $\alpha\in[0,2\pi)$. One can also see that, if we denote $\lambda=\lambda(\bF)$ the supremum of the support of $\bF$ then
\[
 \mu(\bF)\leq \bE(\tau_e)\nu <\lambda\nu\,
\]
(one must assume that $\bF$ is not concentrate in one point, which is the case since $\bF$ is continuous).

We shall also use the following lemma, which is (5.2) of Lemma 5.2 in \cite{hn01}:
\begin{lemma}\label{5.2}
For $\xi\in (0,1)$ and $r>0$ let $A_{\xi,r}$ be the event that there exists $\x\in\bR^2$ with $|\x|\leq 2r$ and $|\x-\v(\x)|\geq r^\xi$. Then, for some constant $c_1>0$,
\[
\bP\big(A_{\xi,r}\big)\leq c_1 e^{-r^{2\xi}}
\]
\end{lemma}

\subsection{$\delta$-straightness of semi-infinite paths}\label{pre-del}
Recall that for $\alpha\in[0,2\pi)$ we have defined that a self-avoiding path $(\x_{n})_{n\geq 1}$, with vertices in $\bR^2$ and such that $|\x_n|\to\infty$ when $n\to\infty$, is a $\alpha$-path if
\[
\lim_{n\to\infty}\frac{\x_{n}}{|\x_{n}|}=e^{i\alpha}:=(\cos\alpha,\sin\alpha)\,, 
\]
and that a sufficient condition for a path $(\x_{n})_{n\geq 1}$ to be a $\alpha$-path for some $\alpha\in[0,2\pi)$ is that, for some fixed $\delta\in(0,1)$ and $c>0$, and for large enough $n$ 
\[
 ang(\x_{n},\x_m)\leq |\x_n|^{-\delta}\,\mbox{ whenever }m>n
\]
($\delta$-straightness). A sufficient condition for $\delta$-straightness is given by the next lemma, which is exactly Lemma 2.7 in \cite{hn01}:
\begin{lemma}\label{hn-101}
If $(\x_n)_\bN$ is a sequence of points in $\bR^2$ with $|\x_n|\to\infty$ when $n\to\infty$, and such that for all large $n$
\[
|\x_{n+1}-\x_n|\leq |\x_n|^{1-\delta}\mbox{ and }d(\x_n,[\x_1,\x_m])\leq |\x_m|^{1-\delta}\mbox{ for }m>n\,,
\]
then there exist a contant $c>0$ such that for all $n$ sufficiently large,
\begin{equation}\label{eang}
ang(\x_n,\x_m)\leq c|\x_n|^{-\delta}\,\mbox{ whenever }m>n\,.
\end{equation}
\end{lemma}

We also consider the $\delta$-straightness property for trees (we have the tree of infection in mind) as follows. For $\epsilon\in[0,\pi)$ let
\[
\C(\x,\epsilon):=\{\y\in\bR^{2}\backslash\{0\}\,:\,ang(\y,\x)\leq\epsilon\}\,. 
\]
If $\calT$ is a tree embedded in $\bR^{2}$, for each
pair $\v,\tilde{\v}\in\calT$ let $\calR_{out}(v,\tilde{v})$ be the set of all
$\hat{\v}\in\calT$ such that the unique path in $\calT$ connecting $\v$ to
$\hat{\v}$ touches $\tilde{\v}$. For $\delta\in(0,1)$, define that $\calT$ is $\delta$-straight at $\v$ if, for
all but finitely many $\tilde{\v}\in\calT$, 
\[
\calR_{out}(\v,\tilde{\v})\subseteq \:\v+C(\tilde{\v}-\v,c|\tilde{\v}-\v|^{-\delta})\,. 
\]
We say that a subset $\calP$ of $\bR^2$ is omnidirectional if, for all $M>0$, the set composed of unit vectors $\v/|\v|$ with $\v\in\calP$ and $|\v|>M$ is dense in $\S^1$. The above lemma, which is Proposition 2.8 in  \cite{hn01}, states that $\delta$-straightness implies existence of an asymptotic orientation:
\begin{lemma}\label{hn-201}
Assume that $\calT$ is a tree embedded in $\bR^2$, whose vertex set is locally finite but omnidirectional, and such that every vertex has finite degree. Assume further that for some vertex $\v$, $\calT$ is $\delta$-straight at $\v$. Then $\calT$ satisfies the following:
\begin{enumerate}
\item Every semi-infinite path in $\calT$ starting from $\v$ has an asymptotic orientation;\\
\item For every $\alpha\in[0,2\pi)$ there exist at least one semi-infinite path in $\calT$ starting at $\v$ and with asymptotic orientation $e^{i\alpha}$.
\item Every semi-infinite path $(\v_n)_{n\geq 1}$ in $\calT$ starting from $\v$ is $\delta$-straight about its asymptotic orientation $e^{i\alpha}$, i.e. $ang(\v_n,e^{i\alpha})<c|\v_n|^{-\delta}$ eventually.
\end{enumerate}
\end{lemma}

\section{Semi-infinite geodesics and the Busemann function}\label{geo}

\subsection{Semi-infinite geodesics: existence}
Recall that a path $\rho=(\v_1,\v_2,\dots)$ in $\calD$ is a semi-infinite geodesic if for all $\v_j,\v_k\in\rho$, the path $(\v_j,\v_{j+1},...,\v_k)$ is the unique geodesic connecting $\v_j$ to $\v_k$. Semi-infinite geodesics starting from $\v\in\calD_v$ and with asymptotic orientation $e^{i\alpha}$ are denoted  $\rho_{\v}(\alpha)$. 

\begin{proposition}\label{pexis}
Let $\Omega_2$ be the event that for all semi-infinite geodesic $\rho$ there exists $\alpha=\alpha(\rho )\in [0,2\pi)$ such that $\rho$ is a $\alpha$-path, and that for all $\alpha\in[0,2\pi)$ and for all $\v\in\calD_v$ there exists at least one geodesic starting from $\v$ and with asymptotic orientation $e^{i\alpha}$. Under \reff{a1}, $\bP(\Omega_2)=1$.
\end{proposition}

The first step to show the existence of semi-infinite geodesics and its convergence is the following result on the fluctuations of $T$, which is exactly Corollary 1.1 in  \cite{p-205}:

\begin{lemma}\label{l2}
Under \reff{a1}, for all $\kappa\in(1/2,1)$ there exist constants $\delta,c_j>0$ such that for all $r\geq 1$ and $s\in [c_1 (\log r)^{1/\delta},c_2 r^{\kappa}]$
\[
\bP\big(|T(\0,r\vec{e}_1)-\mu r|\geq sr^{\kappa}\big)\leq e^{-c_3 s^{\delta}}\,.
\]
\end{lemma}

The second step is to parallel Newman and Piza \cite{np95} to prove that the control of the fluctuations of $T$ can give the control of the fluctuations of a minimizing path connecting $\0$ to $r\vec{e}_1$ about the line segment $[\0,r\vec{e}_1]$. Precisely, for $\xi\in (0,1)$ let 
\[
\C_r^\xi:=\{\x\in\bR^2\,:\,d(\x,[\0,r\vec{e}_1])\leq r^\xi\}\,,
\]
where $[\x,\y]$ denotes the line segment connecting $\x$ to $\y$ and $d(\x,\A)$ denotes the euclidean distance between $\x$ and $\A\subseteq\bR^2$. 

\begin{lemma}\label{l1}
For all $\xi\in(3/4,1)$ there exist constants $c,\delta>0$ such that for all $r\geq 1$
\[
\bP\big(\rho(\0,r\vec{e}_1)\not\subseteq\C_r^\xi\big)\leq e^{-cr^\delta}\,.
\]
\end{lemma}

\proofof{Lemma \ref{l1}} Let $\kappa\in(1/2,1)$, $\tilde{\kappa}\in(\kappa,1)$ and set $\xi\:=(\tilde{\kappa}+1)/2$. Let
\[
 \C^{1,\xi}_{r}:=\{\x\in\bR^{2}\backslash \C^{\xi}_{r}\,:\, d(\x,\C^{\xi}_{r})<r^{\xi}\}\,.
\]
Denote by $F_{r}$ the event defined by the following properties:
\begin{itemize}
\item $\v_{\0},\v_{r\vec{e}_1}\in \C_r^{\xi}$
\item for all edges $\e=(\v,\tilde{\v})$ with $|\v|\leq 2r$ or $|\tilde{\v}|\leq 2r$ we have that $|\v-\tilde{\v}|\leq r^{\xi}$.
\end{itemize}
Notice that $F_r^c\subseteq A_{\xi,r/3}$ (as in Lemma \ref{5.2}), and thus
\begin{equation}\label{HN1}
\bP\big(F_r^c\big)\leq c_1e^{-(r/3)^{2\xi}}\,.
\end{equation}

For each $\z\in\bZ^2$ consider the random variable 
\[
 T_{\z}:=\max_{|\v-\z|\leq 1}\{T(\z,\v)\}\,.
\]
We claim that, under \reff{a1}, for some constants $c_2,c_3>0$
\begin{equation}\label{e7t3}
\bP(T_{\z}\geq r^{\kappa})=\bP\big(T_{\0}\geq r^{\kappa}\big)\leq c_2e^{-c_3r^{\kappa}}\,.
\end{equation}
To see this, notice that $T_\0\leq \sum_{\e\in\calE_1}\tau_\e$, where $\calE_1$ is the set of edges $\e=(\v,\bar{\v})$ in $\calD_e$ with $\C_{\v}\cap \B_{\0}^1\neq\emptyset$ or $\C_{\bar{\v}}\cap \B_{\0}^1\neq\emptyset$. By Lemma \ref{lgraph2}, $\bE\big(\exp(a|\calE_1|)\big)< \infty$ for some $a>0$. Combining this with assumption \reff{a1} and the independence between the Poisson point process and the passage time distribution, one obtains \reff{e7t3}.

Now,
\[
 \big[\rho(\0,r\vec{e}_1)\not\subseteq\C_r^\xi\big]\cap F_{r}\subseteq
\]
\begin{equation}\label{e1t3} 
 \big[\exists v\in\calD_v\cap \C_{r}^{1,\xi}\,:\, T(0,\v)+T(\v,r\vec{e}_1)=T(0,r\vec{e}_1)\big]\subseteq A(r)\,,
\end{equation}
where
\[
A(r):=\big[\exists \z\in\bZ^{2}\cap \C_{r}^{1,\xi}\,:\, T(\0,\z)+T(\z,r\vec{e}_1)\leq T(\0,r\vec{e}_1)+2T_{\z}\big]\,.
\]
Let 
\[
 \Delta(\z,r\vec{e}_1):= \mu|\z-r\vec{e}_1|+\mu|\z|-\mu|r\vec{e}_1|\,.
\]
Thus 
\[
 T(\0,\z)+T(\z,r\vec{e}_1)\leq T(\0,r\vec{e}_1)+2T_{\z}
\]
if and only if,
\[
  \Delta(\z,r\vec{e}_1)\leq \big(T(\0,r\vec{e}_1)-\mu r\big)+\big(\mu|\z|-T(\0,\z)\big)+
\]
\[
 \big(\mu|\z-r\vec{e}_1|-T(\z,r\vec{e}_1)\big)+2T_{\z}\,.
\]
This implies that $A(r)\subseteq\cup_{j=0}^{3}A_{j}(r)$, where  
\[
A_{0}(r):=\big[\exists \z\in\bZ^{2}\cap \C_{r}^{1,\xi}\,:\,T_{\z}\geq \frac{\Delta(\z,r\vec{e}_1)}{8}\big]\,,
\]
\[
A_{1}(r):=\big[\exists \z\in\bZ^{2}\cap \C_{r}^{1,\xi}\,:\,|T(\z,r\vec{e}_1)-\mu|\z-r\vec{e}_1||\geq \frac{\Delta(\z,r\vec{e}_1)}{4}\big]\,,
\]
\[
A_{2}(r):=\big[\exists \z\in\bZ^{2}\cap \C_{r}^{1,\xi}\,:\,|T(\0,\z)-\mu|\z||\geq\frac{\Delta(\z,r\vec{e}_1)}{4}\big]\,,
\]
\[
A_{3}(r):=\big[\exists \z\in\bZ^{2}\cap \C_{r}^{1,\xi}\,:\,|T(\0,r\vec{e}_1)-\mu|r\vec{e}_1||\geq\frac{\Delta(\z,r\vec{e}_1)}{4}\big]\,.
\]
Combining this with \reff{e1t3} one gets that 
\begin{equation}\label{e3t3}
 \bP\big(\rho(\0,r\vec{e}_1)\not\subseteq\C_r^\xi\big)\leq \bP\big(F_{r}^{c}\big)+ \sum_{j=0}^{3}\bP\big(A_{j}(r)\big)\,.
\end{equation}

Notice there exist constants $b_{1},b_{2}>0$ such that for sufficiently large $r>0$ and  $\z\in\bZ^{2}\cap \C_{r}^{1,\xi}$ we have that 
\begin{equation}\label{e4*t3}
 b_{1}r^{\tilde{\kappa}}=b_{1}r^{2\xi-1}\leq \Delta(\z,r\vec{e}_1)\leq b_{2}r^{\xi}=b_{2}r^{\frac{\tilde{\kappa}+1}{2}}\,,
\end{equation}
and
\begin{equation}\label{e4t3}
 r^{\xi}\leq |\z|,|\z-r\vec{e}_1|\leq 2r\,.
\end{equation}

Together with Lemma \reff{l2}, \reff{e4*t3} and \reff{e4t3} yield that for some constant $c_1 >0$ 
\begin{equation}\label{e5t3}
 \bP\big(A_{j}(r)\big)\leq e^{-c_{1}r^{\delta}}\,.
\end{equation}  
Combining \reff{e3t3} with \reff{HN1}, \reff{e7t3}, and \reff{e5t3} one can finish the proof of this lemma. \begin{flushright}\endproof\end{flushright}

For $\v\in\calD_v$ let $\calT_{\v}$ be the union over all $\tilde{\v}\in \calD_v$ of the unique geodesic between $\v$ and $\tilde{\v}$ (the tree of infection at $\v$). Therefore, $\calT_{\v}$ is a tree spanning all $\calD_v$. Thus, the third step is to use Lemma \ref{l1} and the concept of \emph{$\delta$-straightness} for trees discussed before.

\proofof{Proposition \ref{pexis}} Combining Lemma \ref{l1} and Lemma \ref{5.2} with the Borel-Cantelli's lemma, one has that for all $\delta=1-\xi\in(0,1/4)$, almost surely, the assumptions of Lemma \ref{hn-101} hold for all semi-infinite path (geodesic) $(\v_n)_{n\geq 1}$ in $\calT_\v$. So, $\calT_\v$ is $\delta$-straight at $\v$. Since, with probability one, a realization of the Poisson point process is omnidirectional, together with Lemma \ref{hn-201} this yields Proposition \ref{pexis}. \begin{flushright}\endproof\end{flushright}

\begin{remark}\label{r-straight}
Let $\xi\in(3/4,1)$. The almost sure $(1-\xi)$-straightness of the tree of infection also implies that for all $\alpha\in[0,2\pi)$, if $(\v_1,\v_2,\dots)$ is a semi-infinite geodesic with asymptotic orientation $e^{i\alpha}$ then 
\[
ang(\v_n,e^{i\alpha})\leq c|\v_n|^{\xi-1}
\]
for sufficiently large $n$.
\end{remark}

\subsection{Semi-infinite geodesics: uniqueness and coalescence} 
Concerning uniqueness of semi-infinite geodesics we have:
\begin{proposition}\label{puni}
For $\alpha\in[0,2\pi)$ let $\Omega_3(\alpha)$ be the event that for all $\v\in \calD_v$ there exists at most one geodesic starting from $\v$ and with asymptotic orientation $e^{i\alpha}$. Assume only that $\bF$ is continuous. Then $\bP\big(\Omega_3(\alpha)\big)=1$
\end{proposition}

\proofof{Proposition \ref{puni}} For $(k,l)\in\bN^2$, let $A_\alpha(k,l)$ be the event that $U_{k,l}\in \calD_v$ (or equivalently, $N_k\geq l$) and there exists two semi-infinite geodesics starting from $\v=U_{k,l}$, with asymptotic orientation $e^{i\alpha}$, and such that after $\v$ they do not intersect each other. Thus, 
\[
\big(\Omega_3(\alpha)\big)^c\subseteq \cup_{(k,l)\in\bN^2}A_\alpha(k,l)\,.
\]
Now, semi-infinite geodesics starting from the same vertex are not allowed to cross each other and, if a semi-infinite geodesics is caught between two semi-infinite geodesics with the same asymptotic orientation $e^{i\alpha}$ then it must have the asymptotic orientation $e^{i\alpha}$ (by planarity). Therefore, if we denote by $d_{\v}$ the degree of the site $\v=U_{k,l}$ then
\[
|\{\alpha\in[0,2\pi)\,:\,\1_{A_\alpha(k,l)}(\omega)=1\}|\leq d_{\v}(\omega)\,.  
\]
($|A|$ is the cardinality of the set $A$). In particular, almost surely,
\[
\int_{[0,2\pi)}\1_{A_\alpha(k,l)}d\alpha =0\,,
\]
and so, by Fubini's theorem,
\[ 
0\leq \int_{[0,2\pi)}\bP\Big(\big(\Omega_3(\alpha)\big)^c\Big)d\alpha=\int_{\Omega}\big(\int_{[0,2\pi)}\1_{\big(\Omega_3(\alpha)\big)^c}d\alpha\big)d\bP\leq 
\]
\[  
\int_{\Omega}\big(\int_{[0,2\pi)}\sum_{(k,l)}\1_{A_\alpha(k,l)}d\alpha\big)d\bP=\int_{\Omega}\big(\sum_{(k,l)}\int_{[0,2\pi)}\1_{A_\alpha(k,l)}d\alpha\big)d\bP=0\,.
\]
Consequently, there exists $I\subseteq[0,2\pi)$ with total Lebesgue measure so that for all $\alpha\in I$, $\bP\big(\Omega_3(\alpha)\big)=1$. Since $\bP\big(\Omega_3(\alpha)\big)$ does not depend on $\alpha$, this yields Proposition \ref{puni}.
\begin{flushright}\endproof\end{flushright} 

The last result we require to construct the Busemann function is the coalescence behavior of semi-infinite geodesics with the same asymptotic direction:

\begin{proposition}\label{pcoal}
For $\alpha\in[0,2\pi)$ let $\Omega_4(\alpha)\subseteq\Omega_3(\alpha)$ be the event that for all $\v,\bar{\v}\in\calD_v$, if  $\rho_{\v}(\alpha)$ and $\rho_{\bar{\v}}(\alpha)$ do exist (and are unique) then they must coalesce, i.e. there exists $\c=\c(\v,\bar{\v},\alpha)\in\calD_v$ such that 
\[
\rho_{\v}(\alpha)=\rho(\v,\c)\cup\rho_{\c}(\alpha)\,\mbox{ and }\,\rho_{\bar{\v}}(\alpha)=\rho(\bar{\v},\c)\cup\rho_{\c}(\alpha)\,. 
\]
Assume only that $\bF$ is continuous. Then $\bP\big(\Omega_4(\alpha)\big)=1$.
\end{proposition}
 
We note that the almost sure statement in Proposition \ref{pcoal} is for fixed $\alpha\in[0,2\pi)$. As we will see later, almost surely, there exists a random direction $\theta$ so that neither uniqueness nor coalescence hold. Indeed, we will show (in part \ref{pr-rand+strai}) that every branch of the competition interface follows one of those random directions for which coalescence does not hold\footnote{For more on the non coalescence of semi-infinite geodesics see Section 1.3 in \cite{hn01}}.

Let $\calS(\alpha)$ denote the union over
all $\v\in\calD_v$ of $\rho_{\v}(\alpha)$. Then $\calS(\alpha)$ is a forest with say
$N(\alpha)$ disjoint trees. Notice that, on $\big[N(\alpha)\leq 1\big]\cap\Omega_3(\alpha)$, there are no site disjoint semi-infinite geodesic with asymptotic orientation $e^{i\alpha}$. So, Proposition \ref{pcoal} will follow if we prove that $\bP\big(N(\alpha)\leq 1\big)=1$. As noted by Licea and Newman \cite{ln96}, in this set up we can apply
the Burton and Keanne \cite{bk89} method. This method requires several steps which we will be organized as independent claims. To state the first one, let $\delta\in\bQ$ (the set of rational numbers) and $\x_{i}=\big(x_i(1),x_i(2)\big),\tilde{\x}_{i}=\big(\tilde{x}_i(1),\tilde{x}_i(2)\big)\in\bQ^{2}$ for $i=1,...,j$ such that $x_1(2)\leq\dots\leq x_j(2)$ and $\tilde{x}_1(2)\leq\dots\leq\tilde{x}_j(2)$. Assume further that $x_i(1)\leq -\delta$ and that $\tilde{x}_i(1)\geq \delta$. Denote by $A_{\delta}(\x_{1},...,\x_{j},\tilde{\x}_{1},...,\tilde{\x}_{j})$ the event determined by the following:
\begin{itemize}
\item at each $D_{\delta}(\x_{i})$ and $D_{\delta}(\tilde{\x}_{i})$ there is an unique vertex $\v_{i}$ and $\tilde{\v}_{i}$ respectively;
\item each $\e_{i}=(v_{i},\tilde{v}_{i})$ is an edge in $\calD_e$ and $\e_{i}\in\rho_{\v_{i}}(0)$;
\item after $\v_{i}$, $\rho_{\v_{i}}(0)$ has vertices only with strictly positive coordinates;
\item all $\rho_{\v_{i}}(0)$ are disjoint.
\end{itemize}

\begin{claim}\label{emodi1}
If 
\[
\bP\big(N(0)\geq 2\big)>0\,
\]
then 
\[
\bP\big(A_{\delta}(\x_{1},\x_2,\x_3,\tilde{\x}_{1},\tilde{\x}_2,\tilde{\x}_{3})\big)>0\,,
\]
for some $\delta\in\bQ$ and $\x_{i},\tilde{\x}_{i}\in\bQ^{2},i=1,2,3$.
\end{claim}

Since $\bQ$ is enumerable, if $0<\bP\big(N(0)\geq 2\big)$ then there exist $\delta\in\bQ$ and
$\x_{1},\x_{2},\tilde{\x}_{1},\tilde{\x}_{2}\in\bQ^{2}$ such that
\[
0<\bP\big(A_{\delta}(\x_{1},\x_{2},\tilde{\x}_{1},\tilde{\x}_{2})\big)\,. 
\]
Let $c_{n}$ be the maximum between the second coordinate of $\x_{2}$ and
$\tilde{\x}_{2}$ and let $c_{s}$ be the minimum between the second coordinate of
$\x_{1}$ and $\tilde{\x}_{1}$. Consider the rectangle 
\[
\R_{0}:= [-\delta,\delta]\times (c_{s}-\delta,c_{n}+\delta)\,. 
\]
Let $\z_{0}$ be the circumcenter of the rectangle $\R_{0}$ and let $M_{0}$ be the
vertical length of $\R_{0}$. For each $l\in\bZ$ set $\z_{l}:= \z_{0}+lM_0(0,1)$. Denote $\R_{l}:= z_{l}+\R_{0}$ and 
\[
A(l):= A_{\delta}(\x^{l}_{1},\x^{l}_{2},\tilde{\x}^{l}_{1},\tilde{\x}^{l}_{2})\,,
\]
where $\x^{l}_{j}:=\x_{j}+\z_{l}\in \R_l$ and $\tilde{\x}^{l}_{j}:=\tilde{\x}_{j}+\z_{l}$. Thus,  $\bP\big(A(l)\big)=\bP\big(A(0)\big)$. By the Fatou's lemma,
\[
 0<\bP\big(A(0)\big)\leq \bP\big(\limsup_{l} A(l)\big)\leq \bP\big(\cup_{l_{1}\ne l_{2}}A(l_{1})\cap A(l_{2})\big)\,.
\]
Therefore, there are $l_{1},l_{2}$ such that 
\[
0<\bP\big(A(l_{1})\cap A(l_{2})\big)\,. 
\]

Without lost of generality assume that $l_1<l_2$. We claim that, in this case, the geodesic starting from $\v_{1}^{l_{1}}$ can not intersect either the geodesic
starting from $\v_{1}^{l_{2}}$ or the geodesic starting from $\v_{2}^{l_{2}}$. This is so 
because otherwise (by planarity) the geodesic
starting at $\v_{1}^{l_{1}}$ would intersect the geodesic starting from
$\v_{2}^{l_{1}}$, which contradicts the definition of $A(l_{1})$. Thus,
\[
A(l_{1})\cap A(l_{2})\subseteq
A_{\delta}(\x_{1}^{l_{1}},\x_{1}^{l_{2}},\x_{2}^{l_{2}},\tilde{\x}_{1}^{l_{1}},\tilde{\x}_{1}^{l_{2}},\tilde{\x}_{2}^{l_{2}})\,
\]
which yields Claim \ref{emodi1}.

The second step is given by the following claim: for $m,k\geq 0$ let $F_{m,k}$ be the event that some tree in $\calS(0)$ touches a vertex in the rectangle
\[
\R_{m,k}:=\big\{(x(1),x(2))\,:\,0\leq x(1)\leq m\mbox{ and }|x(2)|\leq k\big\}, 
\]
but no other in 
\[
\Q_{m}:=\big\{(x(1),x(2))\,:\,x(1)\leq m\big\}\backslash\R_{m,k}\,.
\]

\begin{claim}\label{lmodi-0}
If for some  $\delta\in\bQ$ and $\x_{i},\tilde{\x}_{i}\in\bQ^{2},i=1,2,3$ we have
\[  
\bP\big(A_{\delta}(\x_{1},\x_2,\x_3,\tilde{\x}_{1},\tilde{\x}_2,\tilde{\x}_{3})\big)>0\,
\]
then 
\[
\bP\big(F_{m,k}\big)>0\,, 
\]
for some $m,k\geq 0$.
\end{claim}

To prove this claim we shall use a local modification argument based on Lemma \ref{lmodif}, and we will divide this proof into two parts: in the first one we will assume that $\bF$ has unbounded support while in the second one we will assume that $\bF$ has bounded support.

\paragraph{\bf Part 1: $\bF$ has unbounded support.} Let $\delta\in\bQ$ and
$\x_{1},\x_{2},\x_{3},\tilde{\x}_{1},\tilde{\x}_{2},\tilde{\x}_{3}\in\bQ^2$ given by Claim \ref{emodi1}. Let
$\R_0:=[-\delta,\delta]\times[c_{s}-\delta,c_{n}+\delta]$, where $c_{n}$ be the maximum between the second coordinate of $\x_{3}$ and
$\tilde{\x}_{3}$ and let $c_{s}$ be the minimum between the second coordinate of
$\x_{1}$ and $\tilde{\x}_{1}$. Denote by $\Xi$ the set of edges which cross the rectangle $\R_0$
and the vertical coordinate axis. Then $\e_{i}:=(\v_{i},\tilde{\v}_{i})\in\Xi$
for all configurations in $A_{\delta}(\x_{1},\x_{2},\x_{3},\tilde{\x}_{1},\tilde{\x}_{2},\tilde{\x}_{3})$ (recall that $\x_{i}\in
C_{\v_{i}}$ and $\tilde{\x}_{i}\in C_{\tilde{\v}_{i}}$). 

Define the event $B_{\lambda}$
by those configurations such that for all $\e=(\v_{1},\v_{2})\in\Xi$ there exists
$\gamma$ with connecting $\v_{1}$ to $\v_{2}$, with $t(\gamma)<\lambda$, but not using edges in  $\Xi$. Since 
\[
\lim_{\lambda\to\infty}\bP\big(B_{\lambda}\big)=1\,, 
\]
we can choose a sufficiently large $\lambda>0$ such that
\begin{equation}\label{emod1}
 \bP\big(A_{\delta}(\x_{1},\x_{2},\x_{3},\tilde{\x}_{1},\tilde{\x}_{2},\tilde{\x}_{3})\cap B_{\lambda}\big)>0\,.
\end{equation}

Now we apply Lemma \ref{lmodif}. To do so define $W(\omega)$, a random element of $\mathnormal{K}$, by the following procedure: given
$\omega\in\Omega$ set
\[
W(\omega):=\big((k_{j},l_{j},n_{j},m_{j})\big)_{j=1,...,q}\, 
\]
by ordering all $(k,l,m,n)$ (according to \reff{E:prescription}) so that $\e(\omega)=\big(U_{k,l}(\omega),U_{m,n}(\omega)\big)\in\Xi(\omega)$ and $\tau_\e\leq\lambda$. Thus $W$ is an ordered representation of the indexes of the edges $\e\in\Xi$ with $\tau_e\leq\lambda$. 

For each $I\in\mathnormal{K}$ let
\[
R_{I}:=(\lambda,+\infty)^{q}\subseteq\Omega_{I}=\bR^q\,,
\]
and let
\[
A:=A_{\delta}(\x_{1},\x_{2},\x_{3},\tilde{\x}_{1},\tilde{\x}_{2},\tilde{\x}_{3})\cap
B_{\lambda}\, 
\]
(given by \reff{emod1}). Since $\bF$ has unbounded support, $\bP_{I}(R_{I})>0$ for all $I\in\mathnormal{K}$. By Lemma \ref{lmodif}, there exist a measurable $\Phi(A)\subseteq\tilde{\Phi}(A)$.

Now consider a configuration $\tilde{\omega}\in\Phi(A)\subseteq\tilde{\Phi}(A)$. By definition, there exists $I\in\mathnormal{K}$, $\omega_1\in\hat{\Omega}_I$, $\omega_2\in \Omega_I$ and $\tilde{\omega}_2\in R_I$ such that $\tilde{\omega}=(\omega_1,\tilde{\omega}_2)$ and $(\omega_1,\omega_2)\in A$. Since $\omega_2$ and $\tilde{\omega}_2$ concern only travel times which are associated to $I$ and $\omega_2\leq\tilde{\omega}_2$ (considering the canonical order in $\bR^q$), the  paths $\rho_{\tilde{\v}_{i}}(0)(\omega_1,\omega_2)$ for $i=1,2,3$ remain disjoint geodesics, with asymptotic orientation $\vec{e}_1$, for the configuration $\tilde{\omega}=(\omega_1,\tilde{\omega}_2)$. By the same reason, $\tilde{\omega}\in B_\lambda$. On the other hand, since $\tilde{\omega}_2\in R_I$, we have that for all $\e\in\Xi$, $\tau_e(\tilde{\omega})>\lambda$ and thus no geodesic could have an edge in $\Xi$. Therefore $\Phi(A)\subseteq F_{m,k}$, where $k:=\max\{c_s,c_n\}$ and $m:=\delta + \max\{\tilde{x}_{1}(1),\tilde{x}_{2}(1),\tilde{x}_{3}(1)\}$. Since  $\bP(A)>0$, we also have that $0<\bP\big(\Phi(A)\big)\leq \bP\big(F_{m,k}\big)$, which yields Claim \ref{lmodi-0} when $\bF$ has unbounded support. 
\paragraph{\bf Part 2: $\bF$ has bounded support.} Consider again $\delta\in\bQ$ and
$\x_{1},\x_{2},\x_{3},\tilde{\x}_{1},\tilde{\x}_{2},\tilde{\x}_{3}\in\bQ^2$ given by Claim \ref{emodi1}. Let  $\vec{e}_2:=(0,1)$, $\c_n:=(0,c_n)$ and $\c_s:=(0,c_s)$. For $\epsilon,\tilde{\epsilon}>0$ and $m>0$, let 
\[
\Q_{m,\tilde{\epsilon}}:=m\vec{e}_1+[-\tilde{\epsilon}m\vec{e}_2,\tilde{\epsilon}\vec{e}_2]\,
\]
and let $B^{\epsilon,\tilde{\epsilon}}_{m}$ be the event that for every $\z\in [\c_s,\c_n]$ and every $\u\in \Q_{m,\tilde{\epsilon}}$,
\begin{equation}\label{ecoal2}
 T(z,u)< (\mu +\epsilon)m\,.
\end{equation}
By the shape theorem, we have that for any $\epsilon>0$ and for sufficiently small $\tilde{\epsilon}$, 
\begin{equation}\label{elimit3}
 \lim_{m\to \infty}\bP\big(B^{\epsilon,\tilde{\epsilon}}_{m}\big)=1\,.
\end{equation}

Denote by $C^{\tilde{\epsilon}}_{m,k}$ the event that for each $i=1,2,3$, $\rho_{\v_{i}}(0)$ touches the hyperplane with direction $\vec{e}_2$ and containing $(0,m)$ for the first time (coming from $\v_{i}$) within the vertical segment $\Q_{m,\tilde{\epsilon}}$. Since all those geodesics are $0$-paths, 
\begin{equation}\label{elimit1}
 \lim_{m\to \infty}\bP\big(C^{\tilde{\epsilon}}_{m}\big)=1
\end{equation}
for all $\tilde{\epsilon}>0$. 

For $m,k>0$ let $C_{m,k}$ denote the event that for each $i=1,2,3$, $\rho_{\v_{i}}(0)$ does not intersect the region consisting of points $(x(1),x(2))\in\bR^2$ such that $x(1)\in [0,m]$ and $|x(2)|>k$. Thus, for any fixed $m>0$, 
\begin{equation}\label{elimit2}
 \lim_{k\to \infty}\bP\big(C_{m,k}\big)=1\,
\end{equation}
(by the same reason to obtain \reff{elimit1}).

Let $\x,\y\in\bR^2$ and recall the definition of the path $\gamma(\x,\y)$ given in Section \ref{pre} (part \ref{pre-geom}). By Lemma \ref{lgraph2},
\begin{equation}\label{egraph2}
\lim_{n\to\infty}\bP\big(|\gamma(\0,n\vec{e}_1)|\geq c_1 n\big)=0\,,
\end{equation} 
for some contant $c_1>0$. We also have considered the graph metric $T_\calD$ and, by Lemma \ref{lgraph1}, 
\begin{equation}\label{egraph1}    
\lim_{n\to\infty}\frac{T_\calD\big([\c_s,\c_n],\Q_{m,\tilde{\epsilon}}\big)}{m}=\nu\,.
\end{equation}

For each $i=1,2,3$, let $\rho_{i}$ denote the piece of $\rho_{\v_{i}}(0)$ between $\tilde{\v}_{i}$ and the first time it intersect $[m\vec{e}_1 -\tilde{\epsilon}m\vec{e}_2,m\vec{e}_1+\tilde{\epsilon}m\vec{e}_2]$, say at the point $\u_i$. For $\z\in[\c_s,\c_n]$ and $\u\in\Q_{m,\tilde{\epsilon}}$, let $\phi(\z,\u)$ be the path connecting $\z$ to $\u$, which first moves vertically by using $\gamma(\z,\v_1)$, then follows $\rho_{1}$, then moves vertically again by using $\gamma(\u_1,\u)$. Thus, on the intersection between $A_{\delta}(\x_{1},\x_{2},\x_{3},\tilde{\x}_{1},\tilde{\x}_{2},\tilde{\x}_{3})$, $C^{\tilde{\epsilon}}_{m}$, $ C_{m,k}$ and $B^{\epsilon,\tilde{\epsilon}}_{m}$, we have that 
\[
 t\big(\phi(z,u)\big)=t\big(\gamma(\z,\v_1)\big)+t\big(\rho_{1}\big)+t\big(\gamma(\u_1,\u)\big)\leq 
\]
\begin{equation}\label{egraph3}
\lambda |\gamma(\c_s,\c_n)|+(\mu+\epsilon)m+\lambda|\gamma(m\vec{e}_1 -\tilde{\epsilon}m\vec{e}_2,m\vec{e}_1+\tilde{\epsilon}m\vec{e}_2 )|\,. 
\end{equation}

We also have that, by \reff{egraph2} and \reff{egraph1} (since $\mu(\bF)<\lambda(\bF)\nu$), there exists $\epsilon_{0},\tilde{\epsilon}_{0}>0$ such that for all $\epsilon<\epsilon_{0}$, $\tilde{\epsilon}<\tilde{\epsilon}_{0}$,
\begin{equation}\label{metric}
 \lim_{m\to \infty}\bP\big(D(\lambda,\epsilon,\tilde{\epsilon})\big)=1\,.
\end{equation}
where $D(\lambda,\epsilon,\tilde{\epsilon})$ is the event that 
\[
\lambda |\gamma(\c_s,\c_n)|+(\mu+\epsilon)m+\lambda|\gamma(m\vec{e}_1 -\tilde{\epsilon}m\vec{e}_2,m\vec{e}_1+\tilde{\epsilon}m\vec{e}_2 )|
\]
\[
\leq (\lambda-\epsilon)T_\calD\big([\c_s,\c_n],\Q_{m,\tilde{\epsilon}}\big)\,.
\]

Let 
\[
A:=A_{\delta}(\x_{1},\x_{2},\x_{3},\tilde{\x}_{1},\tilde{\x}_{2},\tilde{\x}_{3})\cap C^{\tilde{\epsilon}}_{m}\cap C_{m,k}\cap B^{\epsilon,\tilde{\epsilon}}_{m}\cap D(\lambda,\epsilon,\tilde{\epsilon})\,.
\]
Combining \reff{elimit3} with \reff{elimit1}, \reff{elimit2} and \reff{metric}, we get that $\bP\big(A\big)>0$ for sufficiently small $\epsilon>0$ and $\tilde{\epsilon}>0$ and for sufficiently large $m>0$ and $k>0$. Notice that for all configurations in $A$, and every $\z\in[\c_s,\c_n]$ and $\u\in\Q_{m,\tilde{\epsilon}}$ we must have that 
\begin{equation}\label{emod2}
T(\z,\u)\leq t\big(\phi(\z,\u)\big)\leq (\lambda-\epsilon) T_\calD\big([\c_s,\c_n],\Q_{m,\tilde{\epsilon}}\big)\,. 
\end{equation}

Now we are able to use Lemma \ref{lmodif} again. Let $\Xi$ be the set of edges in the interior of the region bounded by
$\rho_{1}$, $\rho_{3}$,  $[\c_s,\c_n]$ and $\Q_{m,\tilde{\epsilon}}$. Define
$W(\omega)$ as follows: given $\omega\in\Omega$ we set
\[
W(\omega):=\big((k_{j},l_{j},m_{j},n_{j})\big)_{j=1,...,q}
\]
by ordering all $(k,l,m,n)$ (according to \reff{E:prescription}) so that
$\e(\omega)=\big(U_{k,l}(\omega),U_{m,n}(\omega)\big)\in\Xi(\omega)$  with $\tau_\e\leq\lambda-\epsilon$. So $W$ represents the indexes of the edges
$\e\in\Xi$ with $\tau_\e\leq\lambda -\epsilon$. For each $I\in\mathnormal{K}$, let
$R_{I}:=(\lambda-\epsilon,\lambda)^{q}\subseteq\Omega_{I}$ and take $A$ above defined. Since $\bF(\lambda-\epsilon)<1$ then $\bP_{I}(R_{I})>0$. Thus, by Lemma \reff{lmodif} there exists a measurable $\Phi(A)\subseteq\tilde{\Phi}(A)$.

Pick a  configuration $\tilde{\omega}=(\omega_1,\omega_2)\in\tilde{\Phi}(A)$. By using the same argument we have done for the other case, one can see that the  paths $\rho_{\tilde{\v}_{i}}(0)(\omega_1,\omega_2)$ for $i=1,3$ remain disjoint geodesics, with asymptotic orientation $\vec{e}_1$, for the configuration $\tilde{\omega}$. The same holds for $\rho_{\u_2}(0)$ and for the inequality \reff{emod2}. On the other hand, by \reff{emod2}, no path $\rho$ connecting $\z\in[\c_s,\c_n]$ to $\u\in\Q_{m,\tilde{\epsilon}}$ that is entirely containing in the region $\Xi$ can be a geodesic for the configuration $\tilde{\omega}$ because, otherwise, 
\[
T(\z,\u)=t(\rho)>(\lambda-\epsilon) T_\calD\big([\c_s,\c_n],\Q_{m,\tilde{\epsilon}}\big)\,. 
\]
This allows us to conclude that
\[
\Phi(A)\subseteq\tilde{\Phi}(A)\subseteq F_{m,k}\,
\]
(with $m,k>0$ given by the definition of $A$). Since $\bP\big(A\big)>0$ we have that $0<\bP\big(\Phi(A)\big)<\bP\big(F_{m,k}\big)$, which yields Claim \ref{lmodi-0} when $\bF$ has bounded support.

The third and last step is: 

\begin{claim}\label{lmodi-1}
$\bP\big(F_{m,k}\big)=0$ for all $m,k\geq 0$.
\end{claim} 

In fact, consider a rectangular array of non-intersecting translates $\R_{m,k}^{\z}$ of the basic rectangle $\R_{m,k}=\R^{\0}_{m,k}$ and of $\Q_m=\Q_m^{\0}$ indexed by $\z\in\bZ^{2}$, and also consider the corresponding event $F_{m,k}^{\z}$. Notice that if $F_{m,k}^{\z}$ and $F_{m,k}^{\tilde{\z}}$ occur, then the corresponding trees in $\calS(0)$ must be disjoint. Thus, if $N_{L}$ is the number of $\z\in [0,L]^{2}$ such that $F_{m,k}^{\z}$ occurs, then
\[
 N_{L}\leq |\{\mbox{ edges crossing the boundary of }[0,L]^{2}\}|.
\]
However, by Lemma \ref{lgraph2}, the expected value of the number of edges crossing the boundary of
$[0,L]^{2}$ is of order $L$. By translation invariance,
\[
\bE\big(N_{L}\big)=n_{L}\bP\big(F_{M,k}\big), 
\]
where $n_{L}$ is the number of rectangles $\R_{m,k}^{z}$ intersecting $[0,L]^{2}$. Since $n_L$ is of order $L^2$, the assumption  $\bP\big(F_{m,k}\big)> 0$ leads to a contradiction.

Now we are able to prove Proposition \ref{pcoal}:

\proofof{Proposition \ref{pcoal}} Combining Claim \ref{emodi1} with Claim \ref{lmodi-0} and Claim \ref{lmodi-1} one obtains 
\begin{equation}\label{lcoal} 
\bP\big(N(\alpha)\leq 1\big)=\bP\big(N(0)\leq 1\big)=1\,.
\end{equation}
By noticing that $\Omega_3\cap\big[N(\alpha\big)\leq 1\big]\subseteq \Omega_4(\alpha)$ one can see that Proposition \ref{pcoal} follows from Proposition \ref{puni} together with \reff{lcoal}. 
\begin{flushright}\endproof\end{flushright}

\subsection{Existence and asymptotics for the Busemann function}\label{asyBuse}
The idea to prove Theorem \ref{tBuse-1} is to combine existence, uniqueness and coalescence of semi-infinite geodesics in a fixed direction $e^{i\alpha}$ to show that if $\z_n\to\infty$ along this direction then for sufficiently large $n$ we have 
\[
T(\x,\z_n)-T(\y,\z_n)=T(\x,\c)-T(\y,\c)\,, 
\]
where $\c$ is coalescence point in direction $e^{i\alpha}$ (Proposition \ref{pcoal}). We begin by introducing what we mean by convergence of paths. Assume that $(\gamma_n)_{n\geq 0} $ is a sequence of finite paths with vertices in $\bR^2$, and for each $n\geq 0$ denote $\gamma_n=(\z_0^n,\z_1^n,\dots,\z_{l_n}^n)$. We define that $\gamma_n$ converges to a semi-infinite path $\gamma=(\x_0,\x_1\,\dots)$, and we write $\gamma=\lim_{n\to\infty}\gamma_n$, if for all $k\geq 1$ there exists $n_k\geq 1$ so that $\gamma_n=(\x_0,\x_1,\dots,\x_k,\z_{k+1}^n,\dots,\z_{l_n}^n)$ for all $n\geq n_k$. For each sequence $(\z_n)_{n\geq 0}$ of vertices in $\bR^2$ with $|\z_n|\to\infty$ and $\z\in\bR^2$ we denote $\Pi\big(\z,(\z_n)_{n\geq 0}\big)$ the set of all semi-infinite paths $\rho$ so that there exists a subsequence $(n_j)_{j\geq 0}$ with $\lim_{j\to\infty} \rho(\z,\z_{n_j})=\rho$.

\begin{lemma}\label{geoconv}
Let $\Omega_1$ be the event that, for all $\alpha\in[0,2\pi)$, if $(\z_n)_{n\geq 1}$ has the asymptotic orientation $e^{i\alpha}$ then: i) $\Pi\big(\z,(\z_n)_{n\geq 1}\big)\neq\emptyset$; ii) every $\rho\in\Pi\big(\z,(\z_n)_{n\geq 1}\big)$ is semi-infinite geodesic with the asymptotic orientation $e^{i\alpha}$. Under \reff{a1},  $\bP\big(\Omega_1\big)=1$.
\end{lemma}

\proofof{Lemma \ref{geoconv}} Let $\calT$ be the tree with vertex set $\cup_{n\geq 1} \rho(\z,\z_n)$ and oriented edges $(\u,\v)\in\calD_e$ (in the Delaunay triangulation) so that $\rho(\z,\u)\subseteq\rho(\z,\v)$. Notice that $\calT$ is an infinite tree. Since every vertex in the Delaunay triangulation has finite degree, the same is true for the vertices in $\calT$. Therefore, by a standard compactness argument, $\Pi\big(\z,(\z_n)_{n\geq 1}\big)\neq\emptyset$. To show that every $\rho\in\Pi\big((\z_n)_{n\geq 1}\big)$ has the asymptotic orientation $e^{i\theta}$ consider $\D\subseteq\S^1$ as in the proof of Theorem \ref{t1}. By Proposition \ref{pexis} and Proposition \ref{puni}, almost surely, for all $\beta\in[0,2\pi)$ such that $e^{i\beta}\in \D$ there exists an unique semi-infinite geodesic starting from $\v(\z)$ and with asymptotic orientation $e^{i\beta}$, which we have denoted by $\rho_{\z}(\beta)$. Now, let $\beta_1,\beta_2\in[0,2\pi)$ such that $e^{i\beta_1},e^{i\beta_2}\in \D$. Assume further that, by following the counter-clokwise orientation of $\S^1$, the unit vector $e^{i\alpha}$ is in between the unit vectors $e^{i\beta_1}$ and $e^{i\beta_2}$. Notice that the paths $\rho_{\z}(\beta_1)$ and $\rho_{\z}(\beta_2)$ bifurcate at some point $\v$ and have no further points in common. On the other hand, $(\z_n)_{n\geq 0}$ has the asymptotic orientation $e^{i\alpha}$. Therefore, once $k$ is large enough, $\rho(\z,\z_k)$ should be in between $\rho_{\z}(\beta_1)$ and $\rho_{\z}(\beta_2)$, and thus the same is true for any limit $\rho$. Since $\D$ is dense in $\S^1$, it follows that $\rho$ has the asymptotic orientation $e^{i\alpha}$. \begin{flushright}\endproof\end{flushright}

\proofof{Theorem \ref{tBuse-1}} Consider the intersection between $\Omega_1$ (path convergence, Lemma \ref{geoconv}) and $\Omega_4(\alpha)$ (coalescence and uniqueness of semi-infinite geodesics, Proposition \ref{pcoal}). In this case, if $(\z_n)_{n\geq 1}$ has the asymptotic orientation $e^{i\alpha}$ then $\lim_{n\to\infty}\rho(\x,\z_n)=\rho_\x(\alpha)$. Together with coalescence, this yields that for $\x,\y\in\bR^2$ there exists
$\c=\c(\x,\y,\alpha)\in\calD_v$ and $n_0>0$ such that 
\[
\rho(\x,\z_n)=\rho(\x,\c)\cup\rho(\c,\z_n)\mbox{ and
}\rho(\y,\z_n)=\rho(\y,\c)\cup\rho(\c,\z_n)\,
\]
for all $n\geq n_0$, which implies that 
\[
T(\x,\z_n)-T(\y,\z_n)=T(\x,\c)-T(\y,\c)\,
\]
for all $n\geq n_0$.
\begin{flushright}\endproof\end{flushright}

\proofof{Theorem \ref{tBuse-2}} Let $\bH_{r}^\alpha$ be the hyperplane that pass through $\a_r:=a_r e^{i\alpha}$ and $r\vec{e}_1$, where $a_r=r/\cos\alpha$. Let $\x_{r}$ be the crossing point between the linear interpolation of $\rho_{\0}(\alpha)$ and $\bH_{r}^\alpha$ that maximizes the distance from $\a_r $. We claim that
\begin{equation}\label{E:asycoal0}
 -T(r\vec{e}_{1},\0)\leq H^{\alpha}(r\vec{e}_{1},\0)\leq T(r\vec{e}_{1},\x_{r})-T(\x_{r},\0)\,.
\end{equation}
The left-hand side of \reff{E:asycoal0} follows directly from the triangle inequality for $T$, since  $H^{\alpha}(r\vec{e}_{1},\0)=T(\x,\c_r)-T(\y,\c_r)$ (as in the proof of Theorem \ref{tBuse-1}). To show the right-hand side, notice that if $\x_{r}\not\in\rho(\0,\c_r)$ then $\c_r \in\rho(\0,\x_{r})$ which implies that $\c_r\in\rho(r\vec{e}_{1},\x_{r})$. Thus 
\[
 H^{\alpha}(r\vec{e}_{1},\0)=T(r\vec{e}_{1},\c_r)-T(\0,\c_r)= T(r\vec{e}_{1},\x_{r})-T(\x_{r},\0)\,.
\]
If $\x_{r}\in\rho(\0,\c_r)$ then 
\[
T(\0,\c_r)=T(\0,\x_{r})+T(\x_{r},\c_r)\,. 
\]
Consequently,
\begin{equation}\label{asycoal0*}
H^{\alpha}(r\vec{e}_{1},\0)=T(r\vec{e}_{1},\c_r)-T(\0,\c_r)=\big(T(r\vec{e}_{1},\c_r)-T(\c_r,\x_{r})\big)-T(\0,\x_{r})\,.
\end{equation}
Since (again the triangle inequality) 
\[
T(r\vec{e}_{1},\c_r)-T(\c_r,\x_{r})\leq T(r\vec{e}_{1},\x_{r})\,,
\]
\reff{asycoal0*} yields \reff{E:asycoal0}. 

Now, 
\[
T(r\vec{e}_{1},\x_{r})-T(\x_{r},\0)=
\]
\[
\big(T(r\vec{e}_{1},\x_{r})-\mu|r\vec{e}_{1}-\a_{r}|\big)\Big(\,:=\,I_1(r)\,\Big)
\]
\[
+\big(\mu|\a_r|-T(\x_{r},\0)\big)\,\,\Big(\,:=\,I_2(r)\,\Big)
\]
\[
+\mu|r\vec{e}_{1}-\a_{r}|-\mu|\a_r|\,\,\Big(\,:=\,I_3(r)\,\Big)\,.
\]
By Remark \ref{r-straight}, if we pick $\xi\in(3/4,1)$ then for some constant $c>0$, almost surely, $|\x_r-\a_r|\leq cr^\xi$ for sufficiently large $r$. On the other hand, by the triangle inequality,
\[
|T(\x_r,r\vec{e}_1)-T(\a_r,r\vec{e}_1)|\leq T(\x_r,\a_r)\mbox{ and }|T(\x_r,\0)-T(\a_r,\0)|\leq T(\x_r,\a_r)\,.
\]
Thus
\[
\limsup_{r\to\infty}\frac{|I_1(r)|}{r}\leq \limsup_{r\to\infty}\frac{|T(r\vec{e}_{1},\a_{r})-\mu|r\vec{e}_{1}-\a_{r}||}{r}+\limsup_{r\to\infty}\frac{\max_{|\z-\a_r|\leq cr^\xi}\{T(\a_r,\z)\}}{r}
\]
and
\[
\limsup_{r\to\infty}\frac{|I_2(r)|}{r}\leq \limsup_{r\to\infty}\frac{|T(\0,\a_{r})-\mu|\a_{r}||}{r}+\limsup_{r\to\infty}\frac{\max_{|\z-\a_r|\leq cr^\xi}\{T(\a_r,\z)\}}{r}\,.
\]

Combining Lemma \ref{l2} with translation invariance one gets that for all $\epsilon>0$
\[
\sum_{r\geq 1}\bP\big(|T(r\vec{e}_{1},\a_{r})-\mu|r\vec{e}_{1}-\a_{r}||\geq \epsilon r\big)<\infty\mbox{ and  }\sum_{r\geq 1}\bP\big(|T(\0,\a_{r})-\mu|\a_{r}||\geq \epsilon r\big)<\infty
\]  
Therefore, by Borel-Cantelli's lemma,
\[
\limsup_{r\to\infty}\frac{|T(r\vec{e}_{1},\a_{r})-\mu|r\vec{e}_{1}-\a_{r}||}{r}=0\mbox{ and  }\limsup_{r\to\infty}\frac{|T(\0,\a_{r})-\mu|\a_{r}||}{r}=0\,.
\]
In \cite{p-205} (Lemma 4.3 there) it is proved that, for some constants $c_0,x_0>0$, if $x>x_0$ then  
\[
\bP\big(T(\0,\z)>x|\z|)\leq e^{-c_0 x|\z|}\,.
\]
By noticing that, with high probability, the number of vertices belonging to a ball of radius $cr^{\xi}$ is of order $r^{2\xi}$, one can get that, for all $\epsilon>0$, 
\[
\sum_{r\geq 1}\bP\big(\max_{|\z|\leq cr^\xi}\{T(\0,\z)\}>\epsilon r\big)<\infty\,.
\]
Thus, together with the Borel-Cantelli's lemma (and translation invariance), this yields
\[
\limsup_{r\to\infty}\frac{\max_{|\z-\a_r|\leq cr^\xi}\{T(\a_r,\z)\}}{r}=0\,.
\]
Consequently,
\[
\limsup_{r\to\infty}\frac{|I_1(r)|}{r}=\limsup_r\frac{|I_2(r)|}{r}=0\,.
\]
Since
\[
\lim_{r\to\infty}\frac{I_3(r)}{r}=\mu\frac{\sin\alpha-1}{\cos\alpha}=-\mu\frac{\cos\alpha}{1+\sin\alpha}\,
\]
we finally have that  
\[
\lim_{r\to\infty}\frac{T(r\vec{e}_{1},\x_{r})-T(\x_{r},\0)}{r}=-\mu\frac{\cos\alpha}{1+\sin\alpha}\,.
\]
Together with \reff{E:asycoal0}, this yields Theorem \ref{tBuse-2}. 
\begin{flushright}\endproof\end{flushright}

\subsection{Competition versus coalescence}\label{pr-rand+strai} In this section we give a sketch of the proof of the statements in Remark \ref{rand+strai}. Let $\varphi:=(\z_{1},\z_{2},\dots)$ be a branch of the competition interface. Thus this branch marks the boundary between two different species, say $j_1$ and $j_2$. Assume further that if one moves along $\z_{n},\z_{n+1},\dots$ then on the right hand side we always see species $j_1$ while on the left hand side one see species $j_2$. By Theorem \ref{t1}, this branch has the direction $e^{i\theta}$ for some $\theta=\theta(\varphi)$. For $l=1,2$, let
$(\v_{n}^{l})_{n\geq 1}$ be the sequence of vertices in $\calD_v\cap\B_{\x_{j_l}}$, so that the tile $\C_{\v_{n}^{l}}$ has an edge boundary that belongs to
$\varphi$ 
. Thus, we have that $\v_{n}^{l}$ has the asymptotic orientation $e^{i\theta(\varphi)}$
(since, by Lemma \ref{5.2}, the distance between $\v_n^{l}$ and the corresponding branch of the competition interface is small if compared with $|\v_n|$).
Together with Lemma \ref{geoconv}, this yields that there exists a subsequence $(n_m)_{m\geq 1}$ and a semi-infinite geodesic $\rho_l$, with asymptotic orientation $\theta(\varphi)$, so that $\rho(\x_l,\v_{n_m}^l)\to \rho_l$. Since $\rho(\x_l,\v_n^l)$ is a geodesic connecting two points in $\B_{\x_{j_l}}(\infty)$, we have that $\rho(\x_l,\v_n^l)\subseteq\B_{\x_{j_l}}(\infty)$ and thus $\rho_l\subseteq \B_{\x_{j_l}}(\infty)$.

Consequently, we have two geodesics $\rho_1$ and $\rho_2$ with the same orientation $e^{i\theta(\varphi)}$, but which do not coalesce (because $\rho_i\subseteq\B_{\x_{j_l}}$ for $l=1,2$). By Proposition \ref{pcoal}, this occurs with zero probability which shows the first statement of Remark \ref{rand+strai}.

By Remark \ref{r-straight}, for all $\xi\in(3/4,1)$, $\rho_{1}$ and $\rho_{2}$ are $(1-\xi)$-straight about its asymptotic orientation $e^{i\theta(\varphi)}$. Since $\varphi$ is caught between $\rho_{1}$ and $\rho_{2}$, this also implies that $\varphi$ is $(1-\xi)$-straight about its asymptotic orientation $e^{i\theta(\varphi)}$, which shows the second statement of Remark \ref{rand+strai}.

\paragraph{\bf Acknowledgment}
This work was developed during my doctoral studies \cite{p04} at Impa and I would like to
thank my adviser, Prof. Vladas Sidoravicius, for his dedication and
encouragement during this period. I also thank Prof. Charles Newman for proposing me the problem studied here,  Prof. Thomas Mountford for a careful reading and useful comments about a previous version of this work, and Prof. James Martin for providing me the numerical simulations in Figure \ref{f2}. Finally, I thank the whole administrative staff of IMPA for their assistance and CNPQ for financing my doctoral studies, without which this work would have not been possible.

\end{document}